\documentclass[journal,10pt,letterpaper,twoside,twocolumn]{IEEEtran}

\usepackage{amsmath,amsbsy,amsfonts,amsthm,graphicx,eufrak}
\usepackage[usenames,dvipsnames]{color}
\usepackage{psfrag}
\usepackage{subeqnarray}
\usepackage{arydshln}
\usepackage[noadjust]{cite}
\usepackage{setspace}
\usepackage{hyperref}
\usepackage[normalem]{ulem}
\usepackage{cancel}
\usepackage{bookmark}
\usepackage{enumerate}
\usepackage{pst-grad} 
\usepackage{pst-plot} 
\usepackage{pst-node}
\usepackage{pst-sigsys}

\def\sh{h}

\def\sw{b}
\def\sx{x}
\def\su{u}
\def\sr{r}
\def\sd{d}
\def\se{e}
\def\sy{y}
\def\M{M}

\def\sb{s}
\def\sa{w}

\def\sv{v}

\def\sss{\text{s}}

\def\ssw{\text{\sw}}
\def\sshc{\text{\^\sh}}
\def\ssh{\text{\sh}}
\def\ssw{\text{\sw}}

\def\ssb{\text{\sb}}
\def\ssu{\text{\su}}

\def\opt{\text{opt}}
\def\vICmod{\xi}

\def\ssxw{\text{\sx}_{\text{\sw}}}
\def\bloc{\text{bloc}}
\def\mod{\text{mod}}

\def\vaICopt{{\boldsymbol \psi}_{\opt}}

\def\ssext{\text{e}}

\def\sH{\MakeUppercase {\vh}}
\def\shc{\ensuremath{\hat{\sh}}}
\def\syc{\ensuremath{\hat{\sy}}}
\def\shc{\ensuremath{\hat{\sh}}}

\def\fCal{{\boldsymbol f}}

\def\NAEC{\ensuremath{N_{\text{AEC}}}}
\def\NBF{\ensuremath{N_{\text{BF}}}}
\def\Nh{\ensuremath{N_{\ssh}}}
\def\Jmin{\ensuremath{J_{\min}}}
\def\Jex{\ensuremath{J_{\text{ex}}}}

\def\F{\ensuremath{F}}
\def\NfCal{\ensuremath{N_{f}}}
\def\NvaIC{\ensuremath{N_{\psi}}}

\def\mBM{{\boldsymbol B}}
\def\vaICw{{\boldsymbol \psi}_b}
\def\vaIChhat{{\boldsymbol \psi}_{\shc}}
\def\vvIC{{\boldsymbol \vartheta}}
\def\vvICmod{{\boldsymbol \vICmod}}
\def\mRbloc{{\boldsymbol R}_{\bloc}}
\def\mRmod{{\boldsymbol R}_{\mod}}
\def\vlambda{\boldsymbol \lambda}
\def\vsigma{\boldsymbol \sigma}
\def\vk{\boldsymbol \nu}

\def\vu {{\boldsymbol \su}}
\def\vx {{\boldsymbol \sx}}
\def\vr{{\boldsymbol \sr}}
\def\vh{{\boldsymbol \sh}}

\def\vhc{\ensuremath {\boldsymbol{ \shc}}}
\def\mH{{\boldsymbol \sH}}

\def\mB{{\boldsymbol {\cal H}}}
\def\mI{{\boldsymbol I}}
\def\mO{{\boldsymbol 0}}
\def\vb{{\boldsymbol \sb}}
\def\va{{\boldsymbol \sa}}

\def\vaIC{{\boldsymbol \psi}}
\def\aopt{{\boldsymbol \sa}_{\opt}}

\def\vw{\boldsymbol b} 
\def\mRbb{{\boldsymbol R}_{\ssb\ssb}}
\def\miRbb{\ensuremath{{\boldsymbol R}^{-1}_{\ssb\ssb}}}

\def\mC{{\boldsymbol C}}

\def\mM{{\boldsymbol {\cal M}}}
\def\mL{{\boldsymbol L}}

\def\mRxx{{\boldsymbol R}_{\ssxw\ssxw}}

\def\vf{{\boldsymbol q}}

\def\vv{{\boldsymbol \sv}}

\def\ftr{{\text {tr}}}

\def\mK{{\boldsymbol {\mathfrak R}}_{\vICmod \vICmod}}
\def\ve{{\boldsymbol \se}}

\def\mRuhcuhc{{\boldsymbol R}_{\ssu_\sshc \ssu_\sshc}}

\def\mCe{{\boldsymbol C}_\ssext}
\def\mBe{{\boldsymbol B}_\ssext}

\def\vqe{{\boldsymbol q}_{\ssext}}

\def\mRvICvIC{{\boldsymbol R}_{\vartheta \vartheta}}

\def\vuh{\vu_{\ssh}}
\def\vuhc{\vu_{\sshc}}
\def\vxw {\vx_{\ssw}}
\def\vxs {\vx_{\sss}}
\def\ves {\ve_{\sss}}
\def\vrs {\vr_{\sss}}
\providecommand{\vwsm}[1]{\vw_{\sss_{#1}}}
\def\vrw {\vr_{\ssw}}

\def\muAEC { \mu_{\mathrm{AEC}}}
\def\muBF { \mu_{\mathrm{BF}}}

\def\mQ{{\boldsymbol Q}}

\def\mLambda{{\boldsymbol \Lambda}}

\def\mRvICvICmod{{\boldsymbol R}_{\vICmod \vICmod}}

\def\vphi{{\boldsymbol \alpha}}

\def\slambdaBk{\lambda_{\Mtransest_k}}

\def\Mtransest{{\boldsymbol \Phi}}

\def\vzero{\boldsymbol 0}

\usepackage[caption=false,font=footnotesize]{subfig}

\begin{document}
%
\title{Statistical Analysis of a GSC-based Jointly Optimized Beamformer-Assisted Acoustic Echo Canceler}
%
%
%


\author{Marcos~H.~Maruo,~\IEEEmembership{Student Member,~IEEE,} Jos\'{e} C. M. Bermudez,~\IEEEmembership{Senior Member,~IEEE} and~Leonardo~S.~Resende,~\IEEEmembership{Member,~IEEE}%
\thanks{This work was supported by CNPq under grant No 140640/2009-6}%
\thanks{This work was partly supported by CNPq under grants No 305377/2009-4 and 473123/2009-6}}%

\maketitle

\begin{abstract}
This work presents a statistical analysis of a class of jointly optimized beamformer-assisted acoustic echo cancelers (AEC) with the beamformer (BF) implemented in the Generalized Sidelobe Canceler (GSC) form and using the least-mean square (LMS) algorithm. The analysis considers the possibility of independent convergence control for the BF and the AEC. The resulting models permit the study of system performance under typical handling of double-talk and channel changes. We show that the joint optimization of the BF-AEC  is equivalent to a linearly-constrained minimum variance problem. Hence, the derived analytical model can be used to predict the transient performance of general adaptive wideband beamformers. We study the transient and steady-state behaviors of the residual mean echo power for stationary Gaussian inputs. A convergence analysis leads to stability bounds for the step-size matrix and design guidelines are derived from the analytical models.   Monte Carlo simulations illustrate the accuracy of the theoretical models and the applicability of the proposed design guidelines. Examples include operation under mild degrees of nonstationarity. 
Finally, we show how a high convergence rate can be achieved using a quasi-Newton  adaptation scheme in which the step-size matrix is designed to whiten the combined input vector.
\end{abstract}


%
\IEEEpeerreviewmaketitle

\section{Introduction}

\IEEEPARstart{A}{coustic} echoes arise in hands free communications when a microphone picks up both the signal radiated in a direct path by a loudspeaker and its reflections at the borders of a reverberant environment. Acoustic echoes tend to degrade intelligibility and listening comfort~\cite{Breining:1999:AECAVHOF,Hansler:2004:AENC}. Modern solutions incorporate adaptive echo cancellers. However, typical room reverberation times require adaptive acoustic echo cancelers with very long responses~\cite{Hansler:2004:AENC,Breining:1999:AECAVHOF}. Also, signal contamination by speech from other talkers, noise and their reflections in the acoustic environment make it difficult to obtain fast convergence and satisfactory echo cancellation with such long cancelers~\cite{Widrow:1976:SNLCLMSAF, Breining:1999:AECAVHOF,Manolakis:2000:SAS, Haykin:1991:AFT, Hansler:2004:AENC}.  
Moreover, conventional acoustic echo cancellation also requires a complex control logic to avoid divergence during double-talk periods~\cite{Bershad:2006:ECLRTDTCC, Tourneret:2009:ECTGLRTDTCC}. Very few studies consider the adaptation during those periods. A recent work~\cite{Gunther:2012:LEPDCDTUSBSS} proposes the use of blind source separation techniques.  Though promising, such technique still lacks computationally efficient solutions.


Assuming it is possible to estimate the direction of arrival (DOA) of the desired speaker, spatial filtering (beamforming) can help attenuate interfering signals in other directions than the desired one. 
Beamformers (BFs) have limited echo suppression capacity due to limits in the array directivity~\cite{VanVeen:1988:BAVASF} and the large number of microphones necessary to suppress all reflections outside the desired DOA~\cite{VanTrees:2002:OAP}.

Acoustic echo cancellation solutions in which BFs and acoustic echo cancelers (AECs) have complementary functions have raised a lot of interest recently~\cite{Kellermann:1997:SCAECABMA,Herbordtt:2000:GSAEC,Herbordt:2001:LGSCEAEC,Herbodt:2004:JOLCMVBAEC,Herbordtt:2005:JOLCMVBAECASR,Kammeyer:2005:NACECB,Maruo:2013:SAJOAECBF-AECS,Maruo:2013:SBAGKLMSA,Guo:2011:AAFECMMSLSUPTFM, Guo:2011:AAFECAGMMSLS, Guo:2011:CMMSLAFECS}.  
BFs and AECs contribute by different means to reduce the residual echo. Hence, using both techniques in a synergistic way can improve the acoustic echo cancellation performance~\cite{Maruo:2011:OTOSBAAEC,Kammeyer:2005:NACECB,Maruo:2013:SAJOAECBF-AECS,Maruo:2013:SBAGKLMSA,Guo:2011:AAFECMMSLSUPTFM, Guo:2011:AAFECAGMMSLS, Guo:2011:CMMSLAFECS}.  
BFs and AECs are usually combined by means of two basic structures~\cite{Brandstein:2001:MASPTA,Kellermann:1997:SCAECABMA}.  
The AEC first structure (AEC-BF) employs one AEC per microphone~\cite{Guo:2011:AAFECMMSLSUPTFM,Guo:2011:AAFECAGMMSLS,Guo:2011:CMMSLAFECS}. 
The BF then processes the AEC outputs for spatial filtering.  It requires several long AECs, leading to very high computational costs~\cite{Guo:2011:AAFECMMSLSUPTFM}. 
Moreover, signals outside the desired DOA must be treated as double talk, complicating the design. The BF first (BF-AEC) structure does the spatial filtering first, leaving basically the echo in the desired DOA to be canceled by a single AEC~\cite{Kammeyer:2005:NACECB,Maruo:2013:SAJOAECBF-AECS,Maruo:2013:SBAGKLMSA}.
This structure presents a significantly lower computational complexity when compared to the AEC-BF structure, even considering that the BF impulse response adds to the length of the response to be identified by the AEC~\cite{Kellermann:1997:SCAECABMA}.
However, as a single AEC has to cancel echoes arriving at many microphones and its desired signal is affected by the BF state, the plant identification model is not valid. Therefore, previous theoretical work has to be used carefully when this structure is studied. In addition, since the AEC solution depends on the BF state, an abrupt change in the desired DOA can lead to a degraded performance until the AEC tracks the new solution.

Alternative structures that have been proposed include the use of polynomial approximations in delay-and-sum beamformers~\cite{Hamalainen:2007:AECDSMAS, Myllyla:2008:ABMDSMAS}, the Transfer-Function Generalized Sidelobe Canceler (TF-GSC)~\cite{Reuven:2007:JNRAEC, Reuven:2004:JAECTFGSCFD,Reuven:2007:MAECNRREUTRGSC, Affes:1996:ASSTAMDTS}, AEC sub-modeling~\cite{Burton:2007:ANSCECB},  mutually exclusive adaptation of the BF and AEC~\cite{Beh:2008:CAECABARSIMR} and wave-domain filtering~\cite{Buchner:2004:WDAFAECFDSBWFS,Buchner:2008:AGDWDAFAAEC}.

Optimization of BF-assisted acoustic echo cancellation systems can be based on different performance surfaces, depending on how the BF and the AEC are optimized. One may define the beamformer performance surface from its own local error~\cite{Herbordtt:2000:GSAEC} or use a joint optimization scheme~\cite{Herbodt:2004:JOLCMVBAEC, Maruo:2013:SAJOAECBF-AECS} in which the global cancellation error is used to jointly optimize the BF and the AEC.  The joint optimization scheme was first proposed in~\cite{Herbodt:2004:JOLCMVBAEC}. It was later applied to a robot speech recognition system~\cite{Herbordtt:2005:JOLCMVBAECASR}.  Joint BF-AEC optimization leads to an optimal solution with better echo cancellation performance than separate BF and AEC optimizations~\cite{Maruo:2011:OTOSBAAEC}. 

Despite the possibilities of combined BF and AEC acoustic echo cancellation systems, we find only few analyses of their transient behavior in the literature.  The AEC-BF structure has been studied for the acoustic echo cancellation problem in~\cite{Guo:2011:AAFECMMSLSUPTFM, Guo:2011:AAFECAGMMSLS, Guo:2011:CMMSLAFECS} and for the acoustic feedback cancellation in \cite{Guo:2013:ACLAFCS}.  A stochastic model has been derived using the power transfer function method for the case of a fixed BF, where just the AEC is adapted. More recently, the transient behavior of a system where a direct-form BF and an AEC are jointly adapted using equal and fixed step-sizes was analyzed in~\cite{Maruo:2013:SAJOAECBF-AECS,Maruo:2013:SBAGKLMSA}. 
The derived analytical model was shown to accurately predict the adaptive system behavior and corroborated previous experimental findings that the same cancellation performance of a single-microphone AEC can be achieved with a shorter AEC when the possibility of spatial filtering is available~\cite{Kallinger:SOCMCECWB}. The model, based on the equivalence to a coventional Linearly Constrained Minimum Variance (LCMV) optimization, allows the use of previous analytical results~\cite{Frost:1972:AALCAAP,Godara:1986:ACLMSAWAABUPS}.  

Adaptive LCMV beamforming may be implemented in many different forms and by using different algorithms~\cite{Frost:1972:AALCAAP,Griffiths:1982:AAALCAB,Buckley:1987:SSFLCMVB,Campos:2002:CAAEHT,Buckley:1986:BBGSC,Resende:1996db,Resende:1996us}.
The direct and GSC forms are equivalent in that both lead to the same optimal solution~\cite{Breed:2002:SPELCMVGSCBF}. 
For some algorithms and under specific conditions they are equivalent even in their transient behavior~\cite{Griffiths:1982:AAALCAB,Buckley:1986:BBGSC,Resende:1996db,Resende:1996us,Werner:2003:OERLSILCMVGSC}.  
Both forms tend to have comparable computational complexities for a small number of constraints. However, the GSC form offers greater design flexibility due to the possibility of choosing the block matrix. Good choices may lead reduced computational complexity~\cite[p. 31]{Griffiths:1982:AAALCAB}. Also, robust GSC implementations with an adaptive bock matrix have been proposed to account for small changes in the desired signal DOA~\cite{Hoshuyama:1999:ARABMABMUCAF, Herbodt:2001:CEFDRGSC, ETT:ETT4460130207, Herbordtt:2005:JOLCMVBAECASR, Herbodt:2007:MBWRFDAFIAAB}. Therefore, it is of interest to study the behavior of the GSC form of the BF-AEC structure.

This work extends the analysis in \cite{Maruo:2013:SAJOAECBF-AECS,Maruo:2013:SBAGKLMSA} to the study of the transient behavior of the jointly optimized BF-AEC structure in the GSC form. We formulate the joint optimization as a single constrained optimization problem, what simplifies the statistical analysis.  Moreover, the analysis incorporates the case of a positive-definite step-size matrix~\cite{Mikhael:1986:AFIAP,Rupp:2000:RCLMSATVMSS,Evans:1993:AIVSSAA,Harris:1986:AVSAFA,Dallinger:2009:ASSLAGTA}. The incorporation of this extra flexibility to the model is particularly interesting for BF-assisted echo cancelers, as their AEC adaptation control logic stops AEC adaptation during double-talk periods~\cite{Bershad:2006:ECLRTDTCC, Tourneret:2009:ECTGLRTDTCC}, while the BF continues adapting using a Reference Signal Based (RSB) structure~\cite{Zhang:2011:ABRVCULA} with the AEC output as the reference signal. 
The problem of designing an adaptive filter with step-size matrices was studied in~\cite{Mikhael:1986:AFIAP,Rupp:2000:RCLMSATVMSS,Evans:1993:AIVSSAA,Harris:1986:AVSAFA,Dallinger:2009:ASSLAGTA}. An exponential model for the echo channel and information on the room reverberation time were exploited in~\cite{Makino:1993:EWSNLMSAFBSRIR} to design a step-size optimized algorithm.  In~\cite{Dallinger:2009:ASSLAGTA}, it was shown that LMS algorithm with a step-size matrix is equivalent to the classical LMS algorithm in a transformed space. The same idea is used in our convergence analysis. The analytical model derived in this paper allows the study of the echo canceler behavior including echo-only periods, when AEC adaptation is slower, double-talk periods when only the BF is adapted, and periods after channel changes when fast AEC adaptation is required~\cite{Bershad:2006:ECLRTDTCC, Tourneret:2009:ECTGLRTDTCC}.



The main contributions of this paper are:
\begin{enumerate}[(i)]
 \item  The formulation of the jointly optimized BF-AEC implemented in the GSC form as an LCMV-based GSC. This signal model can be used to design the conventional LCMV-based GSC without loss of generality. Previous theoretical results show that the behavior of the GSC can be studied from the direct form when adaptation uses a single step-size, feasible quiescent solutions and blocking matrices have orthonormal columns~\cite{Griffiths:1982:AAALCAB,Buckley:1986:BBGSC}. Hence the analysis can also be used to design the BF-AEC and conventional LCMV implemented in the direct form using a scalar step-size generalizing the analysis in~\cite{Maruo:2013:SBAGKLMSA}; 
 \item  Incorporation of a step-size matrix. AEC adaptation control logic demands the adaptation of the AEC and BF with different step-sizes during different adaptation scenarios (double-talk, channel changes, tracking, etc)~\cite{Bershad:2006:ECLRTDTCC, Tourneret:2009:ECTGLRTDTCC}. Hence, a novel analysis capable of predicting the transient behavior during different control logic states (different step-sizes) is of undisputable practical relevance. The analysis model uses a positive-definite step-size matrix 
\end{enumerate}
Using the proposed formulation, we derive a statistical model of the behavior of the BF-AEC system implemented in the GSC form with a positive-definite matrix step-size. 
The model also allows the derivation of a high convergence rate algorithm based on a quasi-Newton adaptation scheme in which the step-size matrices are designed to whiten the combined input vector to accelerate convergence.

This paper is organized as follows.  Section~II formulates the problem addressed. Up to Section~II-C the material is basically the same as in~ \cite{Maruo:2013:SBAGKLMSA} and is necessary to establish the notation used in the rest of the paper. Section~II-D introduces the GSC formulation for the problem studied. Section~III describes the analysis structure that allows the analysis of the adaptation using different step-sizes and the quasi-Newton algorithm using the same mathematical framework. Section~IV describes the adaptive solution. Section V derives the statistical model for the adaptive solution. The statistical model convergence is analyzed in Section~VI. Based on the results in section VI, the new quasi-Newton adaptation is derived in Section VII. Section VIII validates the proposed model using simulation examples. Finally, conclusions are presented in Section~IX.  In this paper, plain lowercase or uppercase letters denote scalars, lowercase boldface letters denote column vectors and uppercase boldface letters denote matrices.

\section{Problem Formulation} \label{sec:ProblemFormulation}

\begin{figure}[h!tb]
\begin{center}
 \includegraphics[width=0.45\textwidth]{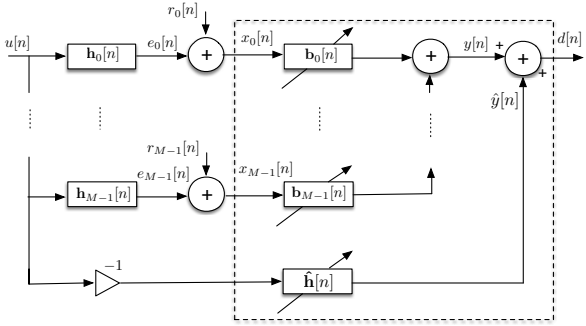}
\caption{BF-AEC system configuration in the direct-form structure~\cite{Maruo:2011:OTOSBAAEC}.}
\label{fig:modelosimplificado}
\end{center}
\end{figure}

Fig.~\ref{fig:modelosimplificado} shows the BF-AEC structure with $\M$ echo impulse response vectors $\vh_m$ of length $\Nh$, $\M$ microphone signals $x_m[n]$, one adaptive wideband beamformer composed of $\M$ filters $\vw_m[n]$ of length $\NBF$ and an adaptive AEC filter $\vhc[n]$ of length $\NAEC$.  
We assume impulse responses $\vh_m$ constant and stationary signals for mathematical tractability~\cite[pp. 348--351]{Haykin:1991:AFT}.  The analysis for a time variant echo path becomes specially challenging in this case even for the simple random walk system nonstationarity model \cite{Haykin:1991:AFT, Manolakis:2000:SAS}. This is because a time variant loudspeaker-enclosure-microphone (LEM) model would lead to a nonstationary beamformer input signal.  Moreover, the statistically independent increments to the channel response vectors $\vh_m$ due to the random walk model would be time-correlated by the BF filters.  This would render the analysis too complex even for such simple nonstationarity model, making it very hard to study fundamental properties of the algorithm behavior.  The study for nonstationary input signals requires a specific model for the input nonstationarity.  To the best of our knowledge, there is no generally accepted model for signal nonstationarity. On the other hand, model predictions derived under stationarity assumptions can still show tendencies of the algorithm behavior for reasonably small degrees of nonstationarity~\cite[p. 595]{Manolakis:2000:SAS}. Simulation results in Section~\ref{subsec:nonstat} will illustrate that this is the case for the present study. 

It has been conjectured that the spatial filtering realized by the BF reduces the required AEC length, as compared to the conventional finite impulse response (FIR) AEC structure~\cite{Kallinger:SOCMCECWB}.  Hence, our analysis considers the possibility of an AEC shorter than the LEM impulse responses by admitting $\NAEC \le \Nh$. 

\subsection{The Beamformer Input Vector}

Each of the $\M$ LEM impulse responses $\vh_m$, $m=0,\ldots,\M-1$, models the transmission of the far-end signal $u[n]$ from the speaker to one of the $\M$ microphones.  The adaptive wideband beamformer is composed by $\M$ FIR filters with impulse responses $\vw_m[n]$,  $m=0,\ldots,\M-1$, each of length $\NBF$~\cite{Liu:2010:WBCT}.  The echo signal at the $m$th microphone is given by~\cite{Breining:1999:AECAVHOF}
\begin{align}
 e_m[n] = \vh_m^T \vuh[n]
 \label{eqn:ekofn}
\end{align}
where
\begin{align}
 \vuh[n] =  \left[ u[n], u[n-1], \cdots, u[n-(\Nh-1)] \right]^T
\end{align}
is the LEM plant input vector.

Grouping the LEM responses as columns of the matrix
\begin{align}
 \mH = \left[ \vh_0 \quad \vh_1 \quad \cdots \quad \vh_{\M-1}  \right]
\label{eqn:matrixH}
\end{align}
and defining the echo snapshot vector as
\begin{align}
\ves[n] = \left[  e_0[n], e_1[n], \cdots, e_{\M-1}[n] \right]^T
\label{eqn:vectores}
\end{align}
\eqref{eqn:ekofn} leads to the linear mapping
\begin{align}
 \ves[n] = \mH^T\vuh[n].
\label{eqn:echosnapshots}
\end{align}

The $m$th microphone signal $x_m[n]$ is the sum of a near-end signal $r_m[n]$ and an echo $e_m[n]$:  
\begin{align}
x_m[n] = e_m[n] + r_m[n], \quad m = 0,\ldots,M-1.
\label{eqn:nearendscalar}
\end{align}
Each signal $r_m[n]$ is composed of local speech, local interferences and random noise.  
 We define the microphone array snapshot $\vxs[n] $ as the vector composed by all $x_m[n]$:
\begin{align}
 \vxs[n] = & \Big[ x_0[n],  x_1[n], \cdots,  x_{\M-1}[n]  \Big]^T.
 \label{eqn:vxs}
\end{align}
Then, combining \eqref{eqn:vectores}, \eqref{eqn:nearendscalar} and \eqref{eqn:vxs}  yields
\begin{align}
 \vxs[n] = \ves[n] + \vrs[n]\nonumber
\end{align}
where 
$\vrs[n] = \left[  r_0[n], r_1[n], \cdots, r_{\M-1}[n] \right]^T$
is the near-end signal component snapshot. 

We now define the extended far-end sample vector as
\begin{align}
 \vu[n] = & \left[ u[n], u[n-1], \cdots, u[n-(\Nh+\NBF-2)] \right]^T
 \label{eqn:extendedfarend}
\end{align}
where the dimension of $\vu[n]$ is the length of the convolution of $\vh_m[n]$ and $\vw_m[n]$. Then to express the microphone array input signals (the echo signals) as functions of $\vu[n]$ we rewrite \eqref{eqn:echosnapshots} as 
\begin{equation}
\ves[n-k]  = \left[ \mO_{\M \times k} \qquad \mH^T \qquad \mO_{\M \times \NBF - (k+1)} \right] \vu[n]
\end{equation}
where $\mO_{N_1 \times N_2}$ denotes the null matrix with dimension $N_1$ lines and $N_2$ columns. 
Then, defining the $\M.\NBF \times 1$ stacked echo vector
\begin{align}
\ve[n]  = \left[ \ves^T[n], \ves^T[n-1], \cdots, \ves^T[n-(\NBF-1)] \right]^T
\label{eqn:BeamformerEcho}
\end{align}
we can write
\begin{align}
 \ve[n]  = \mB^T \vu[n]
 \label{eqn:echovector}
\end{align}
where
\begin{align}
\mB = 
\left[ 
\begin{array}{c:c:c:c}
 				& \mO_{1 \times \M} 		&\cdots 	&  \mO_{\NBF-1 \times M} \\
 \mH 			& \mH 			&\cdots 	& \mH \\
 \mO_{\NBF-1 \times \M} 	& \mO_{\NBF-2 \times \M}	&\cdots & \\
\end{array}
\right]
\label{eqn:matrixB}
\end{align}
is the $\Nh + \NBF-1 \times \M.\NBF$ modified echo channel matrix.  Note that $ \ve[n] $ contains the echo signals for the time window corresponding to the length of the BF impulse response.

Using \eqref{eqn:nearendscalar}, \eqref{eqn:vxs}, \eqref{eqn:BeamformerEcho} and \eqref{eqn:echovector}, and defining the $\M.\NBF \times 1$ near-end vector component (without echo) as
\begin{align}
 \vrw[n] = & \left[ \vrs^T[n], \vrs^T[n-1], \cdots, \vrs^T[n-(\NBF-1)] \right]^T
\end{align}
and the $\M. \NBF \times 1$ combined beamformer input regressor as~\cite{Frost:1972:AALCAAP}
\begin{align}
 \vxw[n] = & \Big[ \vxs^T[n] , \vxs^T[n-1], \cdots, \vxs^T[n-(\NBF-1)]  \Big]^T
\label{eqn:xn}
\end{align}
we write the beamformer input vector as
\begin{align}
 \vxw[n] = \mB^T \vu[n] + \vrw[n].
\label{eqn:nearendvectorform}
\end{align}

\subsection{The Residual Echo}

Define the vector $\vwsm{\ell} [n]$ of the $\ell$th components of all vectors $\vw_{m} [n]$, $m=0,\ldots,\M-1$, at time $n$ as
\begin{align*}
 \vwsm{\ell} [n] = & \Big[ b_{0_\ell}[n] , \cdots, b_{\M-1_\ell}[n] \Big]^T,  \ell=0,\ldots,\NBF-1.
\end{align*}
We then write the beamformer output $y[n]$ as
\begin{align}
y[n] = \sum_{\ell=0}^{\NBF-1} \vxs^T [n-\ell]  \vw_{\sss_\ell} [n].
\end{align}
Now, defining the stacked beamformer weight vector 
\begin{align}
 \vw[n] = & \left[ \vwsm{0}^T[n], \vwsm{1}^T [n], \cdots, \vwsm{\NBF-1}^T[n] \right]^T
\label{eqn:w}
\end{align}
we can write $y[n]$ as the inner product
\begin{align}
y[n]  =  \vxw^T[n] \vw[n].
\label{eq:Expryn}
\end{align}
Next, defining the AEC weight vector
\begin{align}
 \vhc[n]  = & \left[ \shc_0[n], \shc_1[n], \cdots, \shc_{\NAEC-1}[n]  \right]^T
 \label{eqn:hc}
\end{align}
and the AEC input vector
\begin{align}
\vuhc[n] =  \left[ u[n], u[n-1], \cdots, u[n-(\NAEC-1)] \right]^T
 \label{eqn:un}
\end{align}
we can write the AEC output as
\begin{align}
\syc[n] =  \vhc^T[n] \vuhc[n].
\label{eq:Expryhatn}
\end{align}  

Using \eqref{eq:Expryn} and \eqref{eq:Expryhatn} we write the residual echo $\sd[n]$ as the inner product
\begin{align}
\sd[n] = -\vuhc^T[n]\vhc[n]  +\vxw^T[n]\vw[n].
\label{eqn:dninitial}
\end{align}

\section{The Analysis Structure}
With the problem formulation presented in Section~\ref{sec:ProblemFormulation}, we can define an analysis problem that corresponds to the study of a single GSC structure that combines the beamformer and the AEC adaptations. To this end, we define the $\NAEC + \M.\NBF \times 1$ stacked input vector
\begin{align}
\vb[n] = \left[ - \vuhc{^T}[n], \vxw{^T}[n]  \right]^T
\label{eqn:adjointsig}
\end{align}
and, from \eqref{eqn:w} and \eqref{eqn:hc}, the stacked coefficient vector
\begin{align}
\va[n] = \left[ \vhc^T[n], \vw^T[n] \right]^T.
\label{eqn:adjointweightvector}
\end{align}

Then, we can write the residual echo $\sd[n]$ as the inner product
\begin{align}
\sd[n] =\vb^T[n] \va[n].
\label{eqn:dn}
\end{align}


This simple model will permit to relate the study of the BF-AEC structure to that of the LCMV problem.

Interestingly, input vectors $\vu[n]$ in \eqref{eqn:extendedfarend} and $\vuhc[n]$ in \eqref{eqn:un} are related by
\begin{align}
 \vuhc[n] = \left[\begin{array}{cc}\mI_{\NAEC} 	& \mO_{\NAEC \times (\Nh + \NBF - \NAEC -1)} \end{array}\right]\vu[n]
\label{eqn:farendvectorform}
\end{align}
where we use the notation $\mI_{K}$ to denote the $K\times K$ identity matrix. 
Hence, \eqref{eqn:nearendvectorform} and \eqref{eqn:farendvectorform} permit to write $\vb[n]$ in \eqref{eqn:adjointsig} as a function of the input vectors $\vu[n]$ and $\vrw[n]$.  Equation~\eqref{eqn:farendvectorform} also allows to study the algorithm performance for $\Nh > \NAEC$, and thus verifies the possibility of reducing $\NAEC$ by increasing the number of microphones.


\subsection{Performance Surface}
The mean output power (MOP) performance surface $J$ is defined as the mean value of $d^2[n]$ conditioned on $\va[n]=\va$. From \eqref{eqn:dn},
\begin{align}
 J=E\{d^2[n]|\va[n]=\va\} & = E\left\{ \va^T \vb[n] \vb^T[n] \va \right\} \nonumber\\
 & = \va^T \mRbb \va.
\label{eqn:Ed2n}
\end{align}
where $\mRbb = E\{\vb[n]\vb^T[n]\}$ is the input autocorrelation matrix. A set of $\NfCal$ linear constraints on the beamformer coefficients implements the spatial filtering.
Usually, an $\M \NBF \times \NfCal$ constraint matrix  $\mC$ and an $\NfCal \times 1$ response vector $\fCal$ jointly define the 
frequency response 
in the desired
DOA~\cite{Frost:1972:AALCAAP, Buckley:1987:SSFLCMVB}.

To formulate the linear constraints as a function of the combined coefficient vector, we define the extended constraint matrix~\cite{Herbodt:2004:JOLCMVBAEC}
\begin{align}
 \mCe = \left[ \mO_{\NfCal \times \NAEC} \qquad \mC^T \right]^T.
\label{eqn:extendedConstraintMatrix}
\end{align}
Finally, the joint optimization problem can be formulated as
\begin{subequations}
\label{eqn:optproblem}
\begin{equation}
 \aopt = \arg \displaystyle{\min_{\va}}\ \va^T \mRbb \va \label{subeqn:objfunc}
 \end{equation}
 \begin{equation}
 \text{subject to } \mCe^T\va  =  \fCal 
\label{subeqn:constraint}
 \end{equation}
\end{subequations}
and the optimal solution is given by~\cite{Frost:1972:AALCAAP} $\aopt = \miRbb \mCe \left(\mCe ^T \miRbb \mCe \right)^{-1} \fCal$.
\subsection{Implementation using the GSC Form}

\begin{figure}[ht]
\centering
\includegraphics[scale=0.5]{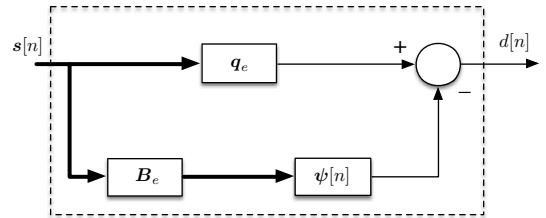}
\caption{BF-AEC system in the GSC configuration.}
\label{fig:diagGSC}
\end{figure}

In the GSC form~\cite{Griffiths:1982:AAALCAB}, the dashed square of Fig.~\ref{fig:modelosimplificado} is replaced by the dashed square of Fig.~\ref{fig:diagGSC}.
Feasible solutions to~\eqref{eqn:optproblem} are decomposed as
~\cite{Griffiths:1982:AAALCAB}
\begin{align}
 \va & = \vqe - \mBe \vaIC
\label{eqn:solutionsequivalence}
\end{align}
where $\vqe$ is any feasible solution to~\eqref{subeqn:constraint},
$\mBe$ is a full column-rank $(\NAEC + \M \NBF) \times \NvaIC$-dimensional blocking matrix orthogonal to $\mCe$ ($\mCe^T\mBe=\vzero$), $\vaIC$ is an $\NvaIC$-dimensional vector and $\NvaIC = \NAEC + \M\NBF - \NfCal$. The minimum norm solution to~\eqref{subeqn:constraint} is 
\begin{align}
\vqe = \mCe (\mCe^T \mCe)^{-1}\fCal. 
\label{eqn:vqe}
\end{align}

\subsubsection{Optimal Solution}
As $\mCe^T\mBe=\vzero$, $\va$ in~\eqref{eqn:solutionsequivalence} satisfies~\eqref{subeqn:constraint} for any $\vaIC$, and~\eqref{eqn:optproblem} becomes an unconstrained optimization problem in $\vaIC$ with solution~\cite{Buckley:1986:BBGSC}
\begin{align}
\vaICopt & = \arg \min_{\vaIC} \vqe^T \mRbb \vqe - 2 \vaIC^T \mBe^T  \mRbb \vqe 
+ \vaIC^T \mRbloc \vaIC
\label{eqn:objfunctionGSC}
\end{align}
where $\mRbloc=\mBe^T \mRbb \mBe$ denotes the blocked input autocorrelation matrix, and from~\eqref{eqn:solutionsequivalence} 
\begin{align}
\aopt = \vqe - \mBe \vaICopt. 
\label{eqn:aoptvaICopt}
\end{align}
Defining the cost function of \eqref{eqn:objfunctionGSC}
\begin{align}
 {\cal C}(\vaIC) = \vqe^T \mRbb \vqe - 2 \vaIC^T \mBe^T  \mRbb \vqe + \vaIC^T \mRbloc \vaIC 
 \label{eqn:costfunctionvaIC}
\end{align}
its gradient with respect to $\vaIC$ is
\begin{align}
 \nabla_{\vaIC}{\cal C}(\vaIC) = - 2 \mBe^T  \mRbb (\vqe - \mBe \vaIC)
\label{eqn:gradobjfunctionGSC}
\end{align}
Setting \eqref{eqn:gradobjfunctionGSC} equal to the null vector yields~\cite{Buckley:1986:BBGSC}
\begin{align}
 \vaICopt = \mRbloc^{-1} \mBe^T  \mRbb \vqe.
\label{eqn:optimalsolGSC}
\end{align}

\section{The Weight Adaptation Equation}

To obtain a model flexible enough to allow the study of the system performance with independent BF and AEC adaptations we choose the following block diagonal form for the blocking matrix
\begin{align}
 \mBe & = \left[ 
\begin{array}{cc}
 -\mI_{\NAEC} & \mO_{\NAEC \times (\M \NBF - \NfCal)} \\
 \mO_{\M \NBF \times \NAEC} & \mBM
\end{array}\right]
\label{eqn:splittedHerbodtBM2}
\end{align}
and split 
\begin{align}
 \vaIC = \left[ \vaIChhat^T, \vaICw^T \right]^T 
 \label{eqn:splitvaIC}
\end{align}
where $[\vaIChhat]_i = [\vaIC]_i$, $i=1,\ldots,\NAEC$ and $[\vaICw]_i = [\vaIC]_{i+\NAEC}$, $i=1,\ldots,\M\NAEC - \NfCal$.
The same block matrix structure has been used in~\cite{Herbodt:2004:JOLCMVBAEC,Herbordtt:2005:JOLCMVBAECASR} for the implementation of the GSC-based BF-AEC acoustic echo canceler.

Using \eqref{eqn:solutionsequivalence} in \eqref{eqn:gradobjfunctionGSC} and noting from \eqref{eqn:dn} that  $\mRbb \va = E\{ \vb[n] \sd[n] \}$ we have
\begin{align}
 \nabla_{\vaIC}{\cal C}(\vaIC) = - 2 \mBe^T  E\{ \vb[n] \sd[n] \}.
\label{eqn:gradobjfunctionGSC2}
\end{align}
Splitting the gradient vector according to~\eqref{eqn:splitvaIC} yields
\begin{align}
 \nabla_{\vaIC}{\cal C}(\vaIC) = \left[ \nabla^T_{\vaIChhat}{\cal C}(\vaIC), \nabla^T_{\vaICw}{\cal C}(\vaIC) \right]^T \nonumber
\end{align}
where, from~\eqref{eqn:adjointsig}, \eqref{eqn:splittedHerbodtBM2}
and \eqref{eqn:gradobjfunctionGSC2}
\begin{subeqnarray}
 \nabla_{\vaIChhat}{\cal C}(\vaIC) = 2 E\{\vuhc[n] \sd[n]\} \slabel{eqn:gradobjfunctionGSChhat}\\
 \nabla_{\vaICw}{\cal C}(\vaIC) = -2 \mBM^T E\{ \vxw[n] \sd[n] \}.\slabel{eqn:gradobjfunctionGSCw}
\end{subeqnarray}

Additionally, setting $\vqe = [\mO_{1 \times \NAEC } \qquad \vf^T]^T$ and 
using~\eqref{eqn:splittedHerbodtBM2} and~\eqref{eqn:vqe} in~\eqref{eqn:solutionsequivalence} yields
\begin{align}
 \va = \left[ \vaIChhat^T, (\vf - \mBM \vaICw)^T \right]^T
 \label{eqn:adjointvectorGSC}
\end{align}
where $\vf = \mC(\mC^T \mC)^{-1}\fCal$.

\begin{figure}[h!tb]
\centering
\includegraphics[width=0.4\textwidth]{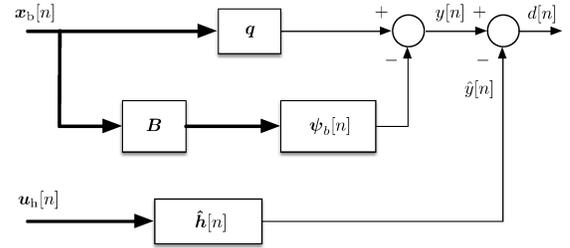}
\caption{Implementation of \eqref{eqn:splitupdates}}
\label{fig:GSC-AEC}
\end{figure}

Comparing \eqref{eqn:adjointweightvector} and \eqref{eqn:adjointvectorGSC} we conclude that
 $\vaIChhat = \vhc$ and $\vw = \vf - \mBM \vaICw$. 
Hence, the steepest-descent algorithms for $\vhc[n]$ and $\vaICw[n]$ with the gradients in \eqref{eqn:gradobjfunctionGSChhat} and~\eqref{eqn:gradobjfunctionGSCw} respectively are
 \begin{subeqnarray}
\label{eqn:splitsteepestdescents} 
  \vhc[n+1] = \vhc[n] - \muAEC E\{\vuhc[n] \sd[n]\} \slabel{eqn:updatehc}\\
  \vaICw[n+1] = \vaICw[n] + \muBF \mBM^T E\{\vxw[n] \sd[n]\} \slabel{eqn:updatevaICw}
 \end{subeqnarray}
where $\muAEC$ and $\muBF$ are the step-size parameters. 
Note that \eqref{eqn:splitsteepestdescents} is different from the steepest descent algorithm for $\vaIC[n]$ unless $\muAEC=\muBF$. This also makes this analysis different from~\cite{Maruo:2013:SBAGKLMSA} by using the equivalence derived on~\cite{Griffiths:1982:AAALCAB,Buckley:1986:BBGSC}. 
However, this extra degree of flexibility is necessary to analyze the behavior of the BF-AEC system under different control logic states that usually act on $(\muAEC, \muBF)$ to avoid divergence.
The stochastic approximations of \eqref{eqn:updatehc} and \eqref{eqn:updatevaICw} yield
\begin{subeqnarray}
\label{eqn:splitupdates}
  \vhc[n+1] \approx \vhc[n] - \muAEC \vuhc[n] \sd[n]\\
  \vaICw[n+1] \approx \vaICw[n] + \muBF \mBM^T \vxw[n] \sd[n].
\end{subeqnarray}
Implementation of \eqref{eqn:splitupdates} has almost the same computational cost of the separate implementation of an LMS implementation of an AEC and a BF demanding only an extra subtraction in the computation of $\sd[n]$. It also requires only one extra memory allocation to account for the second scalar step-size. Despite its simplicity, \eqref{eqn:splitupdates} can model the BF-AEC system behavior under most control logic states. Implementation of \eqref{eqn:splitupdates} is shown in Fig.~\ref{fig:GSC-AEC}.

Finally, the recursive weight update equation is obtained defining the diagonal step-size matrix
\begin{align}
 \mM & = \left[ \begin{array}{cc} \muAEC \mI_{\NAEC} & \mO_{\NAEC \times (\M\NBF - \NfCal)}\\
 \mO_{(\M\NBF -\NfCal) \times \NAEC} & \muBF \mI_{\M\NBF - \NfCal} \end{array} \right]
 \label{eqn:stepsizematrix}
\end{align} 
then \eqref{eqn:splitupdates} can be written as
\begin{align}
 \vaIC[n+1] =  \vaIC[n] + \mM \mBe^T \vb[n] \sd[n].
\label{eqn:updateGSCLMS}
\end{align}
Note that~\eqref{eqn:splitupdates} has the exact same behavior of~\eqref{eqn:updateGSCLMS}, which can be used to study the performance of the practical implementation.

In the following we perform the analysis of an even more general form of \eqref{eqn:updateGSCLMS}, in which the only requirements on $\mM$ and $\mBe$ are that $\mM$ is symmetric positive-definite and $\mBe$ is a full column-rank matrix that satisfies $\mCe^T\mBe=\vzero$.  The typical implementation described above will correspond to a particular case of the more general analysis.

\subsection{Weight Error Vector}

Define the weight error vector $\vv[n] = \va[n] - \aopt$. From~\eqref{eqn:aoptvaICopt},
\begin{equation}
 \vv[n] = \vqe - \mBe \vaIC[n] - \left(\vqe - \mBe \vaICopt\right)  = - \mBe \vvIC[n] 
\label{eqn:mappingweighterrorvectorGSCtodirectform}
\end{equation}
where 
\begin{align}
\vvIC[n]= \vaIC[n] - \vaICopt
\end{align}
denotes the weight error vector of the
unconstrained filter conditioned on
$\mBe$
 and
$\vqe$. From \eqref{eqn:mappingweighterrorvectorGSCtodirectform}, $\vv[n]$ is in the range of $\mBe$.  Hence, $\vv[n]$ is completely determined by $\vvIC[n]$ conditioned on $\mBe$. We then study the behavior of $\vvIC[n]$.


Subtracting $\vaICopt$ from both sides of \eqref{eqn:updateGSCLMS}, using \eqref{eqn:mappingweighterrorvectorGSCtodirectform} with \eqref{eqn:aoptvaICopt} and \eqref{eqn:solutionsequivalence}  
we obtain a recursive update equation for $\vvIC[n]$:
\begin{equation}
\begin{split}
 \vvIC[n+1] = &(\mI_{\NvaIC} - \mM \mBe^T \vb[n] \vb^T[n]\mBe )\vvIC[n] \\
 &+ \mM \mBe^T \vb[n] \vb^T[n]  \aopt.
 \end{split}
\label{eqn:recursionGSC1}
\end{equation}

\section{Statistical Analysis} \label{sec:StatisticalAnalysis}

\subsection{Simplifying Assumptions} 
We now study the behavior of BF-assisted GSC-form echo canceler using \eqref{eqn:updateGSCLMS} under the following typical simplifying assumptions required for mathematical tractability~\cite{Haykin:1991:AFT}
\begin{description}
 \item[A1] $\vb[n]$ is a zero-mean Gaussian vector;
 \item[A2] $\su[n]$ and $\sr[n]$ are statistically independent;
 \item[A3] $\mRbb$ is positive-definite and both $\mCe$ and $\mBe$ have full column rank;    
 \item[A4] The statistical dependence between $\vb[n]\vb^T[n]$ and $\vaIC[n]$ can be neglected;
 \item[A5] The DOA does not change during adaptation.
\end{description}
Though not always valid in practice, these assumptions make analysis viable and frequently lead to results that retain sufficient information to serve as reliable design guidelines~\cite[p. 315]{Haykin:1991:AFT}, \cite{Herbordtt:2000:GSAEC,Herbodt:2004:JOLCMVBAEC}. Simulation results will confirm their reasonability for this analysis. A1 simplifies the evaluation of fourth order moments of $\vb[n]$.  These moments are dependent on the distribution of $\sb[n]$, and the Gaussian distribution combines the advantages of being a good model for several physical processes and simplifying the required mathematical derivations.   A2 is physically reasonable, as $\su[n]$ and $\sr[n]$ are generated at different sides of the communications channel by independent speakers.  A3 is reasonable in practice, as $s[n]$ always has some uncorrelated noise component and both $\mCe$ and $\mBe$ are under reasonable control of the designer.  A4 is required to estimate moments involving the input signal and the weight vector, as the statistical distribution of the latter is unknown.  This assumption is in fact less restrictive than the usually employed independence assumption, which requires $\vb[n]$ and $\vaIC[n]$ to be independent, as discussed in detail in \cite{Minkoff:2001:UASIBRSAFW}. A5 is employed for mathematical tractability and because the main goal of the present analysis is to determine fundamental properties of the adaptive system.

\subsection{Mean Weight Error Vector Behavior}

Taking the expected value of~\eqref{eqn:recursionGSC1} under A4 and using \eqref{eqn:solutionsequivalence} and \eqref{eqn:optimalsolGSC} leads to  

\begin{equation}
 E\{ \vvIC[n+1]\}  =  (\mI_{\NvaIC } - \mM \mRbloc) E\{\vvIC[n]\}
\label{eqn:meanweightvector}
\end{equation}
since
\begin{align}
 E\{\mBe^T \vb[n] \vb^T[n]\} \aopt 
& =  \mBe^T \mRbb \vqe - \mRbloc \mRbloc^{-1}\mBe^T \mRbb \vqe \nonumber\\
& =  \mO_{\NvaIC \times 1}.
\label{eqn:orthogonalityprinc}
\end{align}

Hence, the mean weights converge asymptotically to the optimal solution if all eigenvalues of $\mI_{\NvaIC} - \mM \mRbloc$ are inside the unit circle. In this case,~\eqref{eqn:updateGSCLMS} results in asymptotically unbiased solutions in the mean.

\subsection{Mean Output Power (MOP)}
\label{subsec:secmoment}

To determine the MOP we use \eqref{eqn:dn} with $\va[n] = \vv[n] + \aopt$ and \eqref{eqn:mappingweighterrorvectorGSCtodirectform}. Defining $\mRvICvIC[n]=E\{ \vvIC[n]\vvIC^T[n]\}$ we obtain
\begin{align}
J[n]
 & = E\{ (\aopt - \mBe \vvIC[n])^T \vb[n] \vb^T[n] (\aopt - \mBe \vvIC[n]) \}\nonumber\\
 & = \Jmin + \ftr( \mRvICvIC[n] \mRbloc )
 \label{eqn:Jn}
\end{align}
where we have used  \eqref{eqn:orthogonalityprinc} and A4 to obtain the second line. A recursive expression for the $\NvaIC \times \NvaIC$ matrix $\mRvICvIC[n]$ is derived in the next section to complete the model \eqref{eqn:Jn}. 

\subsection{Correlation Matrix of $\vvIC[n]$}

Post-multiplying \eqref{eqn:recursionGSC1} by its transpose, taking the expected value, using A1--A5  
and 
\eqref{eqn:orthogonalityprinc}
yields
\begin{align}
& \mRvICvIC[n+1]  =  \mRvICvIC[n] 
 - \mM \mRbloc \mRvICvIC[n] - \mRvICvIC[n]\mRbloc \mM & \nonumber\\
 &+ \mM E\Big\{ \mBe^T \vb[n] \vb^T[n] \mBe  \vvIC[n] \vvIC^T[n] \mBe^T \vb[n] \vb^T[n] \mBe \Big\} \mM \nonumber\\
 &+ \Jmin \mM \mRbloc \mM. 
\label{eqn:Rvicvic1}
\end{align}
Using A1, A4 and the Gaussian moment factoring theorem, the expectation in \eqref{eqn:Rvicvic1} is given by
\begin{align}
& E\{ \mBe^T \vb[n] \vb^T[n] \mBe  \vvIC[n] \vvIC^T[n] \mBe^T \vb[n] 
 \vb^T[n] \mBe \} & \nonumber\\
& = 2 \mRbloc \mRvICvIC[n]\mRbloc 
+ \mRbloc \ftr( \mRbloc \mRvICvIC[n] ). &
\label{eqn:termofatmomgauss}
\end{align}
Finally, substituting \eqref{eqn:termofatmomgauss} into
\eqref{eqn:Rvicvic1} yields
\begin{equation}
\begin{split}
\mRvICvIC[n+1]  = &\mRvICvIC[n]   - \mM \mRbloc \mRvICvIC[n] - \mRvICvIC[n] \mRbloc \mM\\
&+ [\Jmin + \ftr( \mRbloc \mRvICvIC[n] )] \mM \mRbloc \mM \\
&+ 2 \mM \mRbloc \mRvICvIC[n]\mRbloc \mM. 
\end{split}
\label{eqn:Rvicvic2}
\end{equation}

Equation~\eqref{eqn:Rvicvic2} completes the MOP model in~\eqref{eqn:Jn}.


\section{Convergence Analysis} \label{sec:ConvergenceAnalysis}

Classical convergence analysis of \eqref{eqn:Rvicvic2} would project $\mRvICvIC[n]$ into the eigenspace of $\mRbloc$ and study the convergence of the diagonal entries of  the transformed matrix~\cite{Manolakis:2000:SAS}.
 The presence of
$\mM$, however, requires a different approach. As $\mM \mRbloc \neq \mRbloc \mM$,  \eqref{eqn:Rvicvic2} is not entirely diagonalizable by the same projection~\cite[p. 558]{Bernstein:2005:MM}. 
Nevertheless, it is still possible to diagonalize both $\mM$ and $\mRbloc$ through contragradient diagonalization~\cite[p. 465]{Bernstein:2005:MM},\cite[p. 466]{Horn:1990:MA}.
As $\mM$ is positive definite, Cholesky decomposition yields $\mM = \mL \mL^T$ with $\mL$ non-singular. Then, we can transform the vector space into
$\vvICmod[n] = \mL^{-1} \vvIC[n]$,
$\mRvICvICmod[n] = E\{\vvICmod[n] \vvICmod^T[n]\} = \mL^{-1} \mRvICvIC[n] \mL^{-T}$ and
\begin{equation}
\mRmod = E\{ \mL^T \mBe^T \vb[n] (\mL^T \mBe^T \vb[n])^T \} = \mL^T \mRbloc \mL. 
\label{eqn:mRmod}
\end{equation}
Hence, pre-multiplying \eqref{eqn:Rvicvic2} by $\mL^{-1}$, post-multiplying by $\mL^{-T}$ and using $ \ftr(\mRvICvICmod[n] \mRmod) =  \ftr(\mRvICvIC[n] \mRbloc)$ yields
\begin{align}
 \mRvICvICmod [n + 1]  =  \mRvICvICmod[n] - \mRmod \mRvICvICmod[n] - \mRvICvICmod[n] \mRmod \nonumber\\
 + \mRmod [\Jmin + \ftr(\mRvICvICmod[n] \mRmod)]  
  + 2 \mRmod \mRvICvICmod[n] \mRmod
\label{eqn:RvICvICmod} 
\end{align}
where $\mRmod$ is symmetric and positive definite. Hence, it is diagonalizable as
$\mRmod = \mQ \mLambda \mQ^T$
with $\mQ^T \mQ = \mI_{\NvaIC}$ and
\begin{align}
 \mLambda = \text{diag} (\lambda_1, \lambda_2, \ldots, \lambda_{\NvaIC}).
\label{eqn:LambdaGSC}
\end{align}
Pre-multiplying \eqref{eqn:RvICvICmod} by $\mQ^T$ and post-multiplying by $\mQ$ yields
\begin{align}
 \mK[n+1] = \mK[n] - \mLambda \mK[n] - \mK[n] \mLambda \nonumber\\
+ \mLambda (\Jmin 
  + \ftr(\mK[n] \mLambda)) + 2 \mLambda \mK[n] \mLambda
\label{eqn:secondorderdiagonalized}
\end{align}
where $\mK[n] = \mQ^T \mRvICvICmod[n] \mQ$.

$\mK[n]$ is an autocorrelation matrix. Then $[\mK[n]]_{i,j}^2 \leq [\mK[n]]_{i,i}[\mK[n]]_{j,j}$, $[\mK[n]]_{i,i} \geq 0$~\cite[p. 251]{Papoulis:1965:PRVSP}, \cite{Prussing:1985:TPMTSDM}, and convergence of~\eqref{eqn:secondorderdiagonalized} can be studied observing only the diagonal elements of $\mK[n]$. Let
$\vk[n]$
be the vector of diagonal entries of $\mK[n]$ and  $\vlambda = [ \lambda_1, \lambda_2, \ldots, \lambda_{\NvaIC}]^T$ be the vector of the eigenvalues of $\mRmod$. Then, from~\eqref{eqn:secondorderdiagonalized}
\begin{equation}
 [\vk[n+1]]_i = \Big[(1 - \lambda_i)^2 + \lambda_i^2\Big][\vk[n]]_i 
+ \lambda_i \big(\vlambda^T \vk[n] + \Jmin\big)
\label{eqn:elementsweightautocorrelation}
\end{equation}
and
\begin{align}
\vk[n+1] = \Mtransest \vk[n] + \Jmin \vlambda
\label{eqn:rhovvrecursiv}
\end{align}
where
$\vlambda^T \vk[n] = \ftr(\mLambda \mK[n])$,
\begin{align}
\Mtransest = \text{diag}(\rho_1, \rho_2, \ldots, \rho_{\NvaIC }) + \vlambda \vlambda^T
\label{eqn:Mtransest}
\end{align}
and $\rho_k = (1 - \lambda_k)^2 + \lambda_k^2$.


The matrix $\Mtransest$ is symmetric and positive definite, as for any nonzero vector $\vphi$ we have
\begin{align}
\vphi^T \Mtransest \vphi = \vphi^T\text{diag}(\rho_1, \rho_2, \ldots, \rho_{\NvaIC})\vphi + \vphi^T \vlambda\vlambda^T \vphi \nonumber\\
 = \sum_{k=1}^{\NvaIC} \left((1 -  \lambda_k)^2 +  \lambda_k^2\right) \left([\vphi]_k\right)^2 + (\vlambda^T \vphi)^2 
> 0. \nonumber
\end{align}

The solution to~\eqref{eqn:rhovvrecursiv} is~\cite{Kailath:1980:LS}
\begin{align}
 \vk[n] & = \Mtransest^n \vk[0] + \Jmin \sum_{j=0}^{n-1} \Mtransest^j \vlambda.
\label{eqn:rhorecursivKailath}
\end{align}
Using \eqref{eqn:rhorecursivKailath} we now study the stability conditions and the steady-state behavior of~\eqref{eqn:updateGSCLMS}.
 
\subsection{Mean Weight Error Revisited}
Pre-multiplying \eqref{eqn:meanweightvector} by $\mL^{-1}$ 
the mean weight error recursion is transformed to~\cite{Rupp:2000:RCLMSATVMSS}
\begin{align}
 E\{ \vvICmod[n+1]\} & = (\mI_{\NvaIC} - \mRmod)E\{ \vvICmod[n]\}.
\label{eqn:meanmodifiedweightvector}
\end{align}
Hence, analysis can be restricted to the eigenspace of $\mRmod$ defined in~\eqref{eqn:mRmod}.

\subsection{Stability Conditions}

Recursion \eqref{eqn:rhorecursivKailath} is a state-space equation whose stability is determined exclusively by the eigenvalues $\slambdaBk$, $k=1,\ldots,\NvaIC$, of $\Mtransest$~\cite{Kailath:1980:LS}. 
From Gershgorin's theorem~\cite{Brualdi:1994:RCP}, 
\begin{align}
 \slambdaBk & < \rho_k + \lambda_k \sum_{j=1}^{\NvaIC} \lambda_j\nonumber\\
&  < 
1 -  2\lambda_k + 2 \lambda_k^2+ \lambda_k \ftr(\mRmod), \quad \forall k
\label{eqn:upperboundeig}
\end{align}
and \eqref{eqn:rhorecursivKailath} is stable if $\slambdaBk<1$ for all $k$. Then, \eqref{eqn:upperboundeig}  leads  to the sufficient condition
\begin{equation}
-2 \lambda_k + 2 \lambda_k^2+ \lambda_k \ftr(\mRmod) < 0, \quad \forall k 
\end{equation}
which implies that $\lambda_k \neq 0$ and
\begin{align}
  2 \max \{\lambda_k\} + \ftr(\mRmod) < 2.
\label{eqn:stabsuffcond1}
\end{align}
In most practical cases, a reliable estimate of the eigenvalues of $\mRmod$ is not available a priori and the upper bound in \eqref{eqn:stabsuffcond1} can not be used. However, using the inequality
\begin{align}
\max \{ \lambda_k\} \leq \ftr(\mRmod) \nonumber
\end{align}
it is possible to derive a tighter upper bound~\cite{Horowitz:1981:PACLMSCNBAA,Feuer:1985:CALMSUGD}
\begin{align}
\ftr(\mRmod) = \sum_{k=1}^{\NvaIC} \lambda_k < \frac{2}{3}.
\label{eqn:stabstuffcondfullmatrix}
\end{align}

In the particularly important implementation using~\eqref{eqn:splitupdates}, $\mBe$ is given by~\eqref{eqn:splittedHerbodtBM2} and $\mM$ by~\eqref{eqn:stepsizematrix}. Hence, if we write 
$\ftr(\mRmod)  = \ftr(\mL^T \mRbloc \mL) = \ftr(\mM \mRbloc)$ 
with the matrices in the partitioned form using 
\eqref{eqn:adjointsig},
\eqref{eqn:splittedHerbodtBM2}, 
\eqref{eqn:stepsizematrix} 
and 
\eqref{eqn:mRmod} 
yields
\begin{align}
\ftr(\mRmod) = \muAEC \ftr(\mRuhcuhc) + \muBF \ftr(\mBM^T \mRxx \mBM) 
\label{eqn:traceRmod}
\end{align}
where 
$\mRuhcuhc = E\{\vuhc[n] \vuhc^T[n]\}$ and $\mRxx = E\{\vxw[n] \vxw^T[n]\}$.
Hence, \eqref{eqn:stabstuffcondfullmatrix} becomes
\begin{align}
\muAEC \ftr(\mRuhcuhc) + \muBF \ftr(\mBM^T \mRxx \mBM) < \frac{2}{3}. 
\label{eqn:stabsuffcond2}
\end{align}

\subsection{Excess MOP}

Using  $\ftr(\mRvICvIC \mRbloc) =  \ftr(  \mRvICvICmod[n]  \mRmod ) = \ftr(  \mK[n]  \mLambda  ) =  \vlambda^T \vk[n]$
in \eqref{eqn:Jn}
we write the MOP as a function of $\vlambda$ and $\vk[n]$
\begin{align}
 J[n] 
& = \Jmin + \vlambda^T\vk[n]. 
\label{eqn:Jmsweighterror}
\end{align}
Thus, the excess MOP is given by~\cite[p. 302]{Haykin:1991:AFT}
\begin{align}
\Jex[n]  =  \vlambda^T\vk[n] 
\label{eqn:Jexdef}
\end{align}

\subsection{Steady-State Excess MOP}
When \eqref{eqn:stabsuffcond1} holds,  $\lim_{n \rightarrow \infty} \vk[n+1] = \vk[n] = \vk[\infty]$ and from 
\eqref{eqn:elementsweightautocorrelation}
and \eqref{eqn:Jexdef} we have
\begin{equation}
 [\vk[\infty]]_{i}  =   \Big[(1 -  \lambda_i)^2 +  \lambda_i^2\Big]  [\vk[\infty]]_{i}  + \lambda_i (\Jex[\infty] + \Jmin) 
\end{equation}
which solved for $[\vk[\infty]]_{i}$ yields
\begin{equation} 
[\vk[\infty]]_{i} = (\Jex[\infty] + \Jmin) \frac{1}{2 - 2 \lambda_i  }.
\end{equation}
Using this result in \eqref{eqn:Jexdef} as $n \rightarrow \infty$ and solving for $\Jex[\infty]$ yields
\begin{equation}
\Jex[\infty]  =  \Jmin \frac{ \frac{1}{2}\sum_{i=1}^{\NvaIC} 
\frac{ \lambda_i}{1 -  \lambda_i  } } { 1 -\frac{1}{2} \sum_{i=1}^{\NvaIC} 
\frac{\lambda_i}{1 -  \lambda_i  } }.
\label{eqn:Jexinftyquasela}
\end{equation}
From 
\eqref{eqn:mRmod}, $\mRbloc=\mBe^T \mRbb \mBe$
and A3, $\mRmod$ is symmetric and positive definite. Hence, its largest eigenvalue is related to its largest singular value through $\max\{\vlambda\} = \sqrt{\max \{\vsigma\}}$ where $\vsigma$ denotes the vector of singular values of $\mRmod$. The largest singular value of a matrix is equal to its $2$-norm~\cite[pg.78]{Chen:1998:LSTD}. Then, using the Cauchy-Schwarz inequality~\cite[pg. 291]{Horn:1990:MA}
\begin{align}
\max \{\vlambda\}  
 & < \left( \|\mL^{T}\|_2 \|\mRbloc\|_2 \|\mL\|_2 \right)^{\frac{1}{2}}\nonumber\\
 & < \left( \|\mM\|_2 \|\mRbloc\|_2 \right)^{\frac{1}{2}}
\end{align}
where both $\mM$ and $\mRbloc$ are symmetric positive definite matrices. 
Hence, for $\max\{\vlambda_{\mM}\} \max\{\vlambda_{R_{\mathrm{bloc}}}\} \ll 1$, where $\vlambda_{\mM}$ and $\vlambda_{R_{\mathrm{bloc}}}$ are vectors containing the eigenvalues of $\mM$ and $\mRbloc$, respectively,
we conclude that $\max\{\vlambda\} \ll 1$ and  
\eqref{eqn:Jexinftyquasela} reduces to
\begin{align}
\Jex[\infty] & \approx \Jmin \frac{\frac{1}{2}\ftr(\mRmod)}{1 - \frac{1}{2} \ftr(\mRmod) }.
\label{eqn:Jexinftyla}
\end{align}

For implementations using~\eqref{eqn:splitupdates}, substituting~\eqref{eqn:traceRmod} in 
\eqref{eqn:Jexinftyla} yields
\begin{align}
\Jex[\infty] & = \Jmin \frac{\muAEC \ftr(\mRuhcuhc) + \muBF \ftr(\mBM^T \mRxx \mBM) }{2 - \muAEC \ftr(\mRuhcuhc) - \muBF \ftr(\mBM^T \mRxx \mBM) }.
\label{eqn:Jexinftyblock}
\end{align}
Further assuming $\muAEC \ftr(\mRuhcuhc) + \muBF \ftr(\mBM^T \mRxx \mBM) \ll 2$ 
we have
\begin{align}
\Jex[\infty] & \approx \frac{\Jmin}{2} \left[\muAEC \ftr(\mRuhcuhc) + \muBF \ftr(\mBM^T \mRxx \mBM)\right] 
\label{eqn:Jexinftysimplificado}
\end{align}

\section{A New Joint Adaptation Algorithm}

The analysis results derived in sections~\ref{sec:StatisticalAnalysis} and \ref{sec:ConvergenceAnalysis} are valid for the general weight update equation \eqref{eqn:updateGSCLMS}.  At the same time, \eqref{eqn:updateGSCLMS} in its general form where $\mM$ and $\mBe$ satisfy only the criteria of $\mM$ being positive-definite and $\mCe^T\mBe=\vzero$ can be considered a new adaptive algorithm that allows weight updating in directions that do not correspond to the stochastic gradient.  



Next, we discuss one possibility of taking advantage of the more flexible structure, namely, designing for a faster convergence speed

The simplest way to guarantee asymptotic convergence to $\vaICopt$ is set 
$\mM=\mu\mI_{\NvaIC}$ in \eqref{eqn:updateGSCLMS} where $\mu=\muAEC=\muBF$, which is the standard LMS update. 
However, it is known that LMS presents a low rate of convergence when the gradient of the performance surface has a low magnitude in the direction of at least one eigenvector of $\mRmod$. 
To alleviate this issue, one may use the step-matrix in~\eqref{eqn:stepsizematrix} or a stochastic approximation of the Newton method. 
The idea underlying quasi-Newton methods is to use an approximation to the inverse Hessian. 
The form of the approximation varies among different methods -- ranging from the simplest, where it remains fixed throughout the iterative process, to the more advanced where improved approximations are built up on the basis of information gathered during the descent process~\cite{Luenberger:1973:ILNP}.

\subsection{High Convergence Rate Block Matrix and Step-Matrix Pair}

Theoretical results show that the rate of convergence of \eqref{eqn:updateGSCLMS} is increased with the reduction eigenvalue spread of $\mL^T\mBe^T \mRbb \mBe \mL$~\cite{Haykin:1991:AFT, Manolakis:2000:SAS, Sayed:2008:AF},  reaching its maximum when all eigenvalues are equal. 
This is because the MOP of an adaptive filter trained with an algorithm of the LMS family decreases over time as a sum of exponentials whose time constants are inversely proportional to the eigenvalues of the autocorrelation matrix of the filter inputs~\cite{Maruo:2013:SBAGKLMSA}. Hence, small eigenvalues create slow convergence modes while large eigenvalues limit the maximum step-size that can be chosen without encountering stability problems as observed in~\eqref{eqn:stabsuffcond1}~\cite{Beaufays:1995:TDAFAA}.  
When all $\mRmod$ eigenvalues are equal, we have
\begin{align}
\mQ^T \mRmod \mQ = \lambda \mI_{\NvaIC}
\label{eqn:hyperconvergence}
\end{align}
where  $\lambda = \ftr(\mRmod)/\NvaIC\approx 2 \frac{\Jex[\infty]}{J[\infty] \NvaIC}$ from~\eqref{eqn:Jexinftyla}. 
Pre-multiplying~\eqref{eqn:hyperconvergence} by $\mQ$, post- multiplying by $\mQ^T$ and noting that $\mQ\mQ^T = \mI_{\NvaIC}$ we conclude that 
$\frac{1}{\lambda} \mRmod =  \mI_{\NvaIC}$. Observing $\mRmod$ structure from~\eqref{eqn:mRmod} yields
 $\mL^T \left( \frac{1}{\lambda} \mRbloc \mL \right) = \mI_{\NvaIC}$ and observing that both $\mL^T$ and $\frac{1}{\lambda}\mRbloc \mL$ are square matrices, we conclude that
 \begin{align}
  \mL^T = \left( \frac{1}{\lambda} \mRbloc \mL \right)^{-1} = \lambda \mL^{-1} \mRbloc^{-1}.
\label{eqn:Lopt}  
 \end{align}
Finally, pre-multiplying~\eqref{eqn:Lopt} by $\mL$ and substituting $\mM = \mL \mL^T$ yields
\begin{gather}
  \mM  = \frac{2}{\NvaIC} \frac{\Jex[\infty]}{J[\infty]} (\mBe^T \mRbb \mBe)^{-1}.  
\label{eqn:optstepsizemat}
\end{gather}
Direct use of \eqref{eqn:optstepsizemat} would require prior knowledge of second-order statistics of $\vb[n]$ and $\Jmin$. Nevertheless, the transient behavior when using~\eqref{eqn:optstepsizemat} is a useful measure of the upper bound on the convergence speed. 
A compromise solution would be the estimation of $\mM$ every few iterations and the use of update~\eqref{eqn:updateGSCLMS}.
In this case, \eqref{eqn:optstepsizemat} becomes a quasi-Newton adaptive filter~\cite{Theodoridis:2001:AFA}.

%
\section{Simulation Examples}
\label{sec:results}

This section presents simulation and design examples to verify the accuracy of the derived model and to illustrate its use in design.  
In all simulations, except when explicitly stated, the far-end signal was drawn from an autoregressive process AR1($a_1$) given by $u[n]=-a_1 u[n-1] + z[n]$, with $z[n]$ a white Gaussian noise with variance $\sigma_z^2$ such that $\sigma_u^2=1$. 
For the accuracy tests, adaptation of the coefficients is assumed to be done during a single-talk period. 
Simulations under different control logic states are considered in subsection~\ref{subsec:controlstates}, where only the BF is adapted during a double-talk period and the convergence of both filters is accelerated after a LEM plant change is detected. 
The LEM plants are designed following the procedure outlined in~\cite{Maruo:2013:SBAGKLMSA} for a uniform linear microphone array oversampled by a factor of $\F=5$\footnote{$\F$ is the ratio between the temporal oversampling factor and the spatial oversampling factor used to generate spatially correlated LEM impulse responses~\cite[Appendix B]{Maruo:2013:SBAGKLMSA}.}, which yields spatially correlated exponential impulse responses. 

\subsection{Model Verification 1}

\begin{figure}[ht]
\subfloat[AR1(-0.9) $f(\muAEC, \muBF) = 2/3$]{
\label{subfig:a}
\centering
\includegraphics[width=0.2\textwidth]{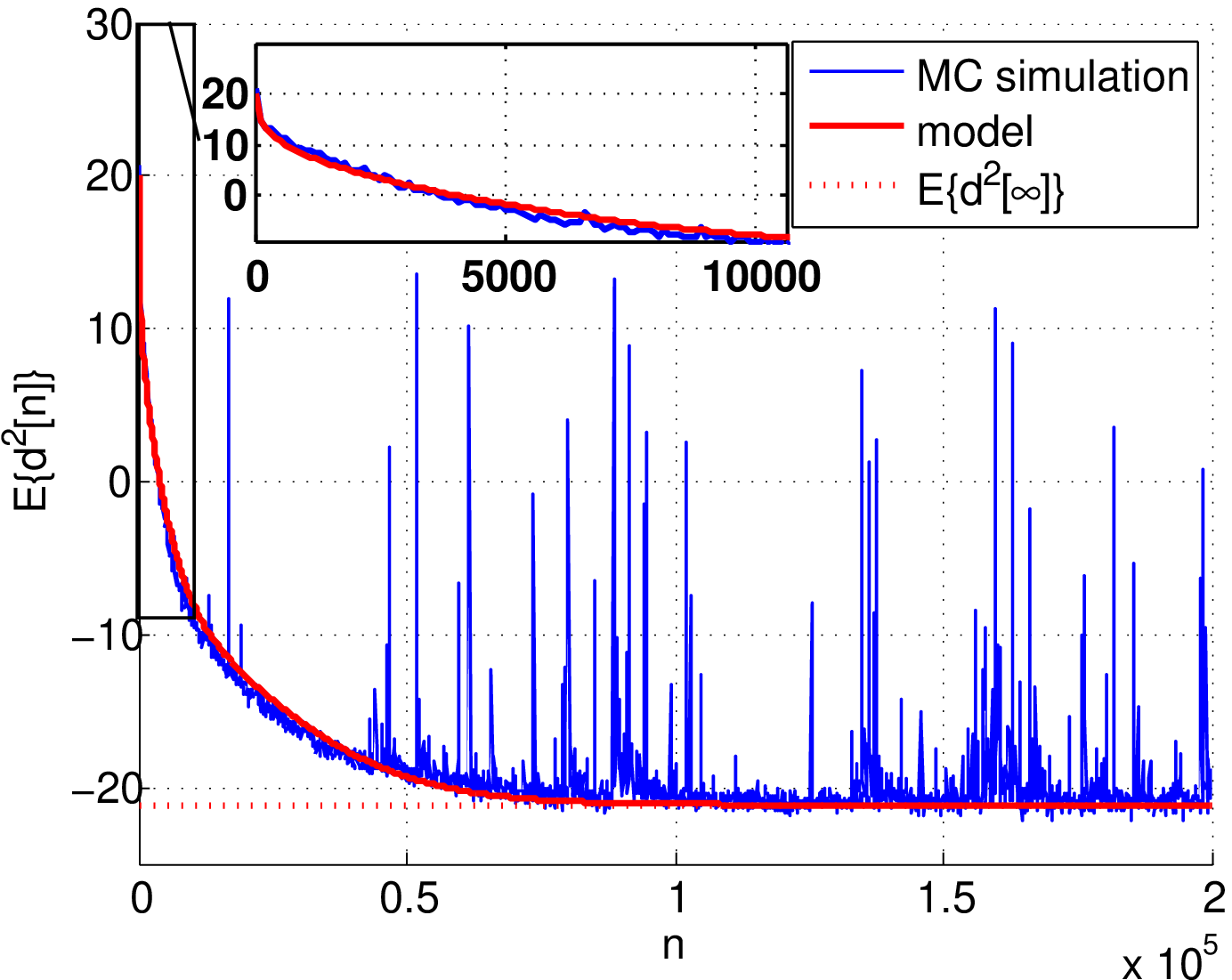}
}
\qquad
\subfloat[AR1(-0.9) $f(\muAEC, \muBF) = 2/30$]{
\label{subfig:b}
\centering
\includegraphics[width=0.2\textwidth]{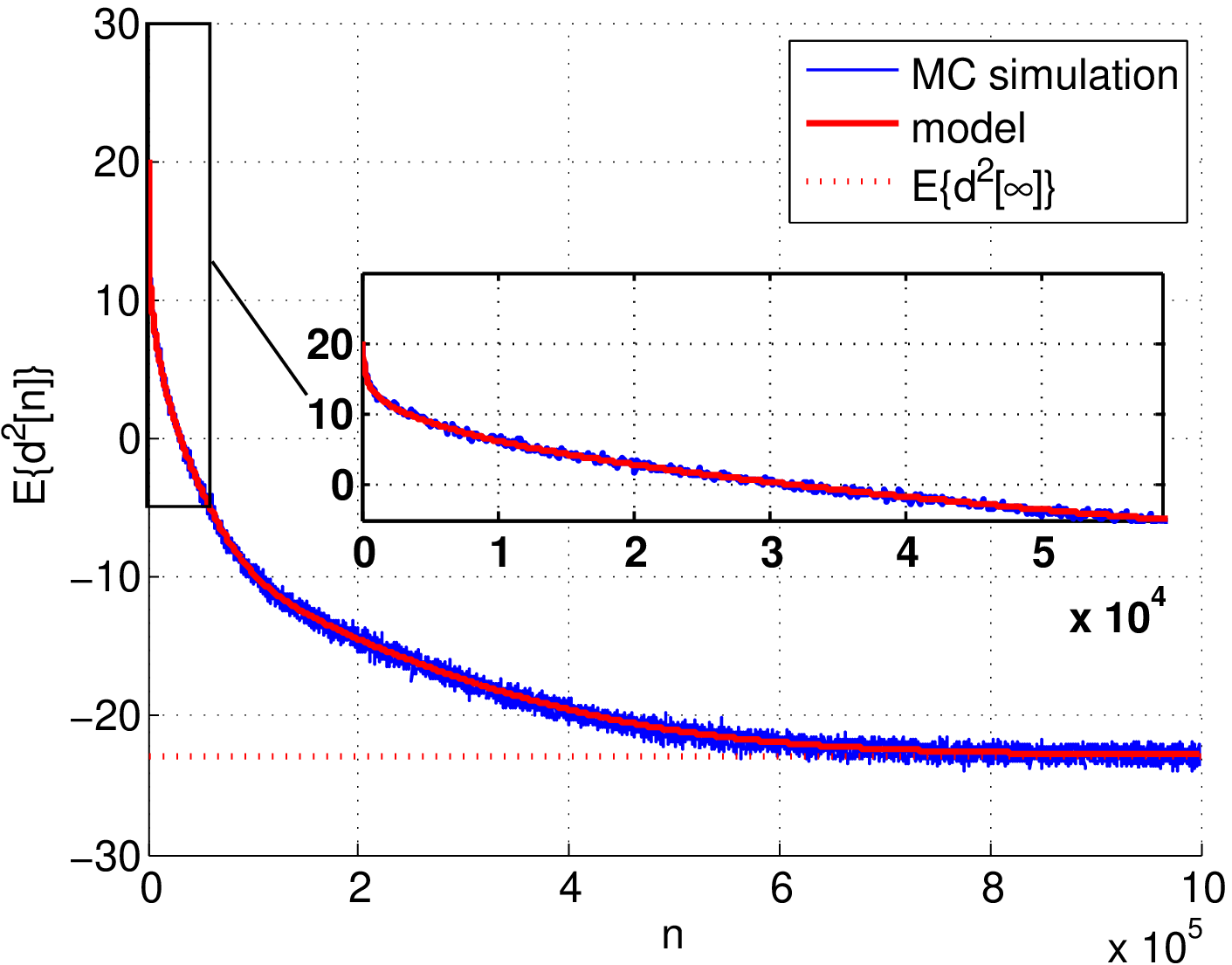}
}\\
\subfloat[AR1(-0.5) $f(\muAEC, \muBF) = 2/3$]{
 \label{subfig:c}
 \centering
\includegraphics[width=0.2\textwidth]{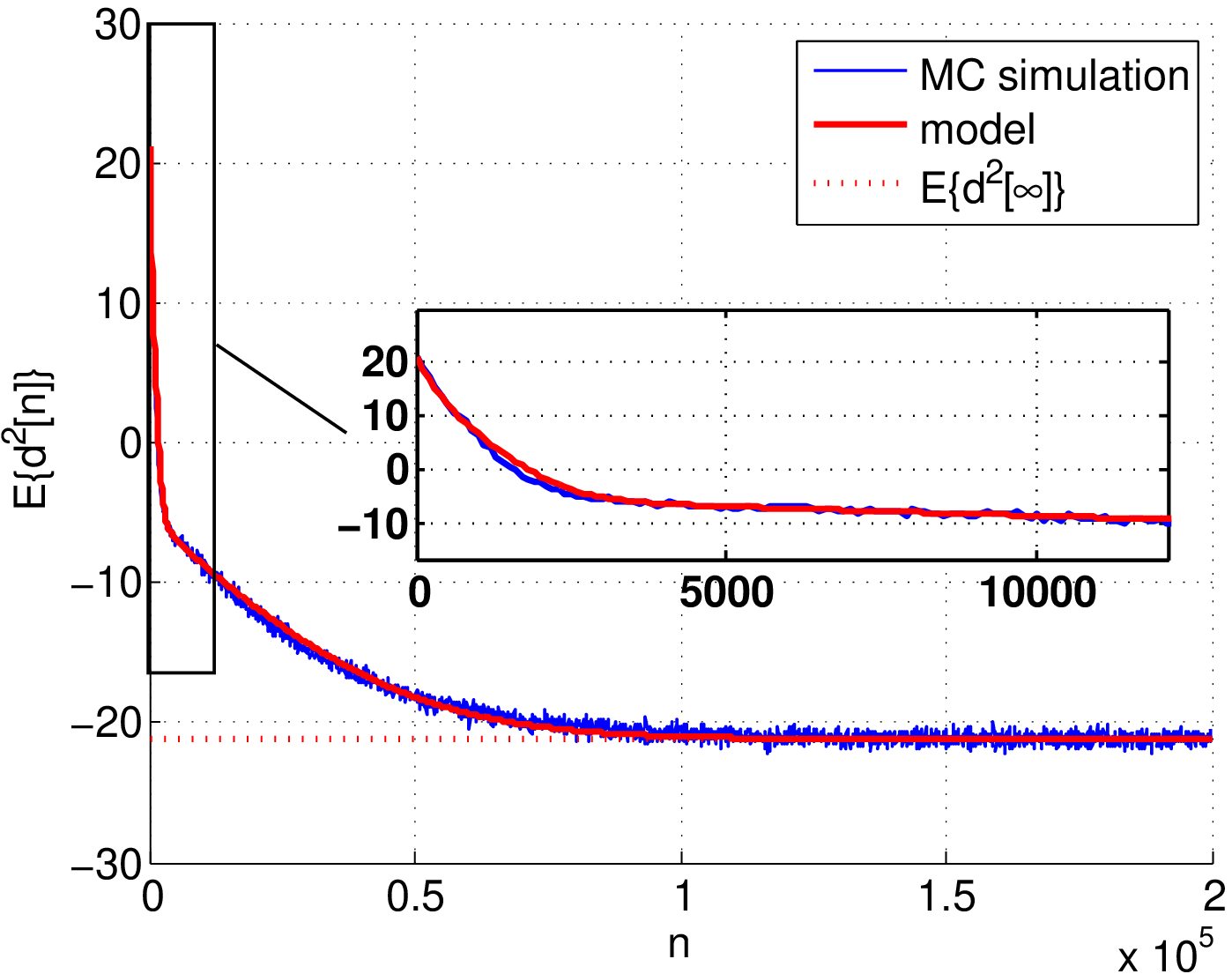}
}
\qquad
\subfloat[AR1(-0.5) $f(\muAEC, \muBF) = 2/30$]{
  \label{subfig:d}
 \centering
\includegraphics[width=0.2\textwidth]{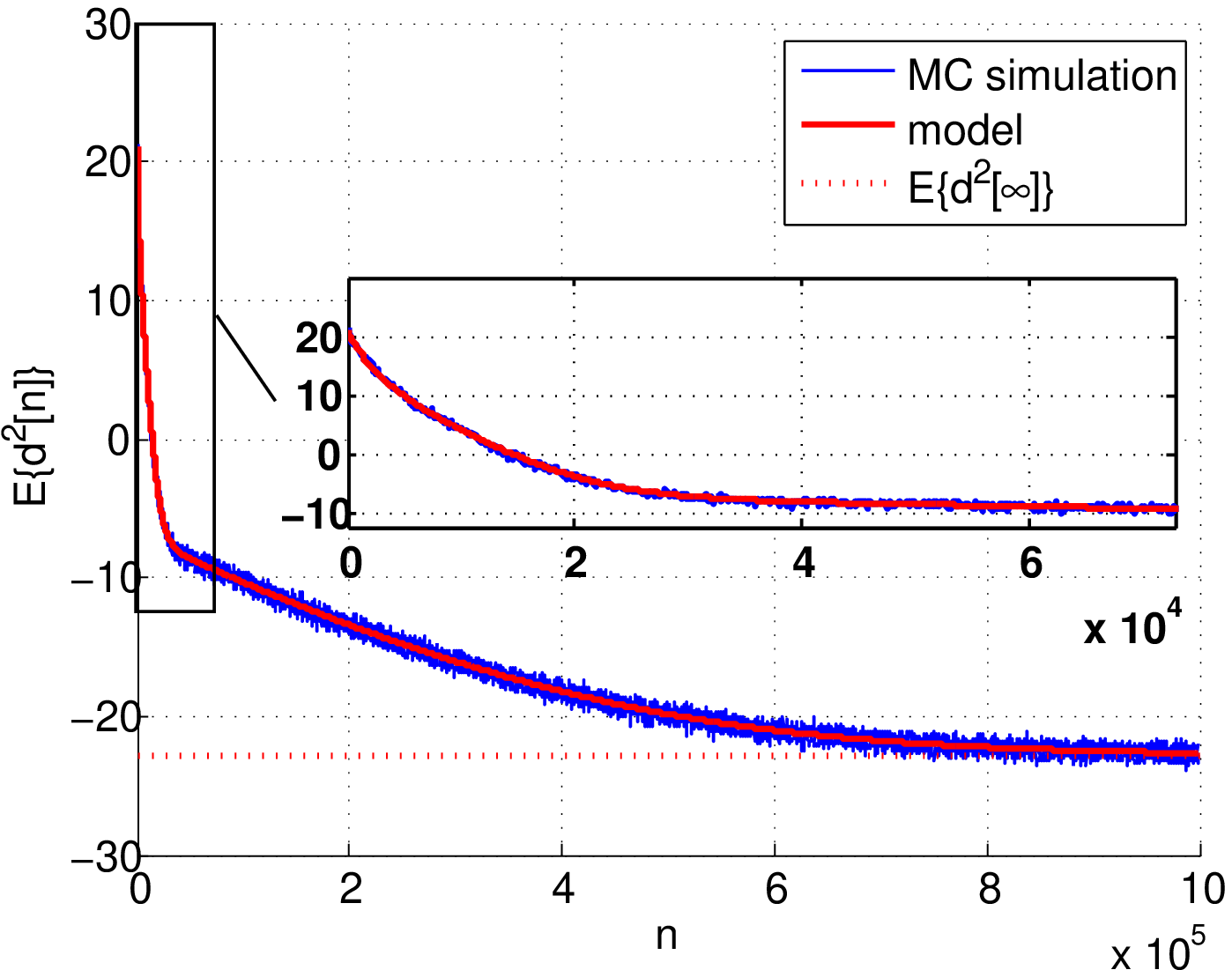}
}\\
%
%
\subfloat[white noise $f(\muAEC, \muBF) = 2/3$]{
\centering
   \includegraphics[width=0.2\textwidth]{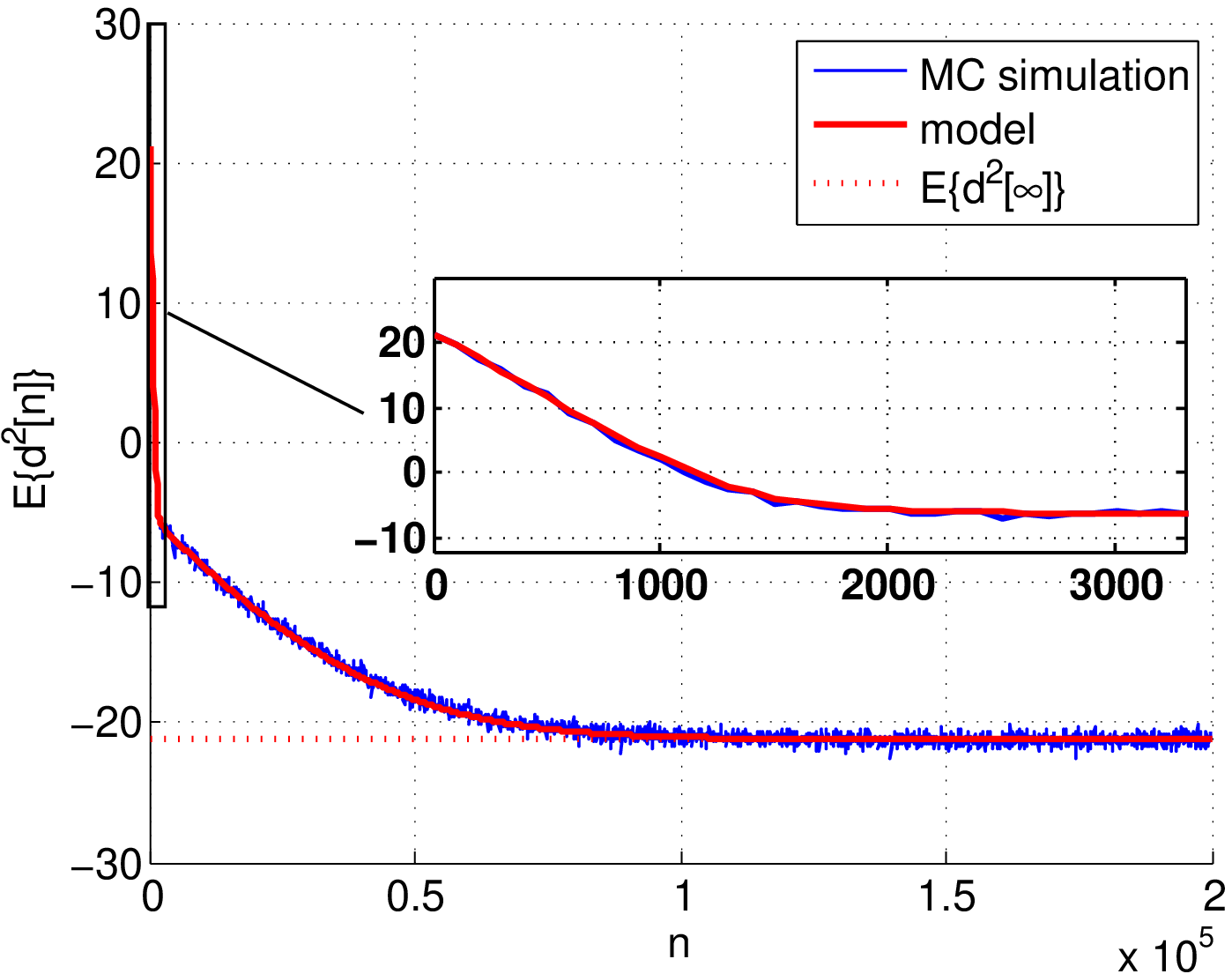}
}
\qquad
\subfloat[white noise $f(\muAEC, \muBF) = 2/30$]{
\centering
   \includegraphics[width=0.2\textwidth]{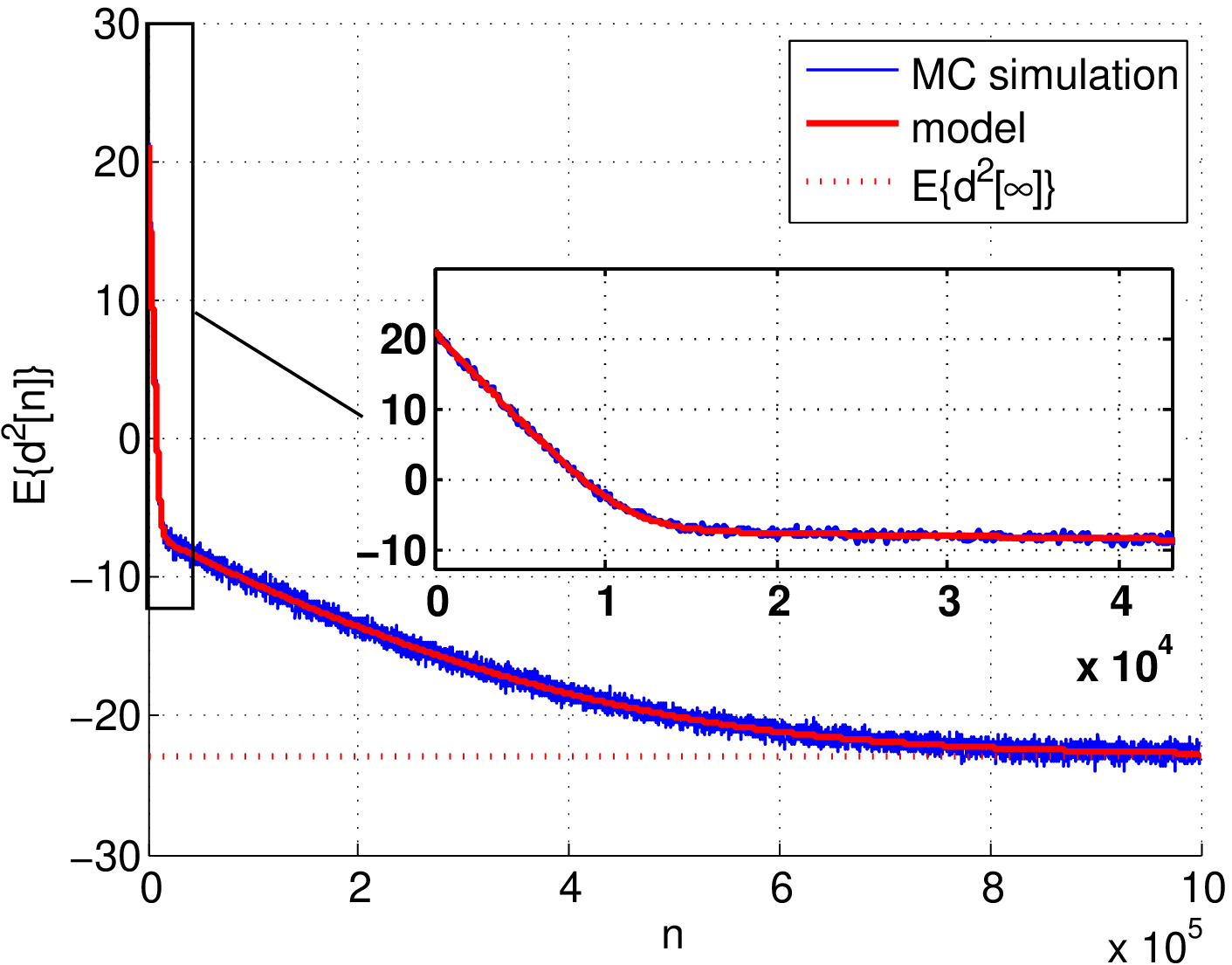}
}
\caption{Proposed model and Monte-Carlo simulation results based on 300 runs for different far-end signal statistics ($\M=2$, $\Nh=128$, $\F=4$, $\NBF=16$ $\NAEC=128$)}
\label{fig:precisionornot}
\end{figure}

The accuracy of the derived model has been verified through Monte-Carlo simulations using several different parameter sets. To conserve space, Fig.~\ref{fig:precisionornot} shows a few of these simulations for $\M=2$ and LEM responses $\vh_1$ and $\vh_2$ with $\Nh=128$. The beamformer filters had $\NBF=16$ and linear phase in the look direction. The AEC length was $\NAEC=128$. 
The noise variance at each microphone was $10^{-2}$.  The theoretical predictions (smooth red curves) are in very good agreement with the Monte Carlo simulations ($300$ runs). Values of $J[\infty]$ using \eqref{eqn:Jexinftyla} are shown by the red horizontal dotted lines.

\subsection{Model Verification 2}

\begin{figure}[htb]
  \centering
\includegraphics[width=8cm]{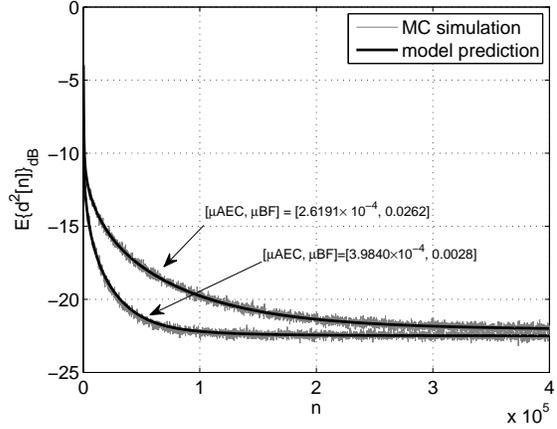}
\caption{Monte-Carlo simulation results for 20 runs with AR1(-0.9) input ($\M=2$, $\Nh=500$, $\NBF=16$ $\NAEC=\Nh + \NBF -1$). }
\label{fig:precisionornot2}
\end{figure}

Consider a unit power first order autorregressive AR1(-0.9) far-end signal
, 2 microphones, $\mathbf{h}_0$ and $\mathbf{h}_1$ with 500 taps each, generated according to the model in~\cite{Maruo:2013:SBAGKLMSA}. The desired DOA was assumed orthogonal to the microphone array. We assumed the absence of double-talk, and noises $r_0[n]$ and $r_1[n]$ were zero-mean white Gaussian with variance $10^{-2}$.  The adaptive BF was designed with $\NBF=16$, linear phase, and all-pass frequency response with $\NfCal = 16$.  The AEC used $\NAEC=\Nh + \NBF -1$.  Fig.~\ref{fig:precisionornot2} shows the predicted and simulated transient MOP.
We tested $2$ scenarios: $[\muAEC, \muBF] = [2.6191\times 10^{-4}, 0.0262]$ and $[\muAEC, \muBF]=[3.9840\times10^{-4}, 0.0028]$. 
Fig.~\ref{fig:precisionornot2} shows excellent agreement between theory and predictions in both cases. Counterintuitively,  the results show that a larger convergence speed does not necessarily imply a higher steady-state error.

\subsection{Model Verification 3}
\label{subsec:nonstat}

\begin{figure}[htb]
\centering
\subfloat[$\eta = 0.01$]{
\centering
\includegraphics[width=0.4\textwidth]{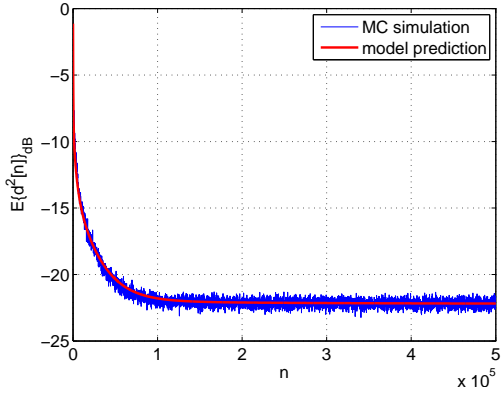}
}\\
\subfloat[$\eta = 0.1$]{
\centering
\includegraphics[width=0.4\textwidth]{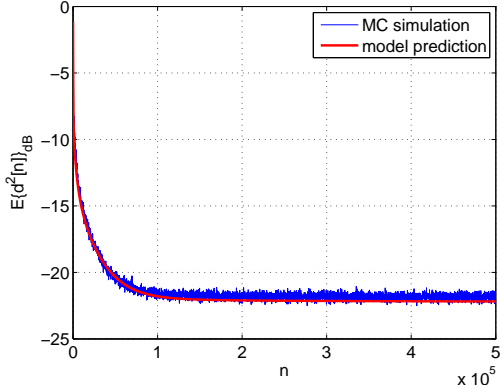}
}\\
\subfloat[$\eta = 0.5$]{
\centering
\includegraphics[width=0.4\textwidth]{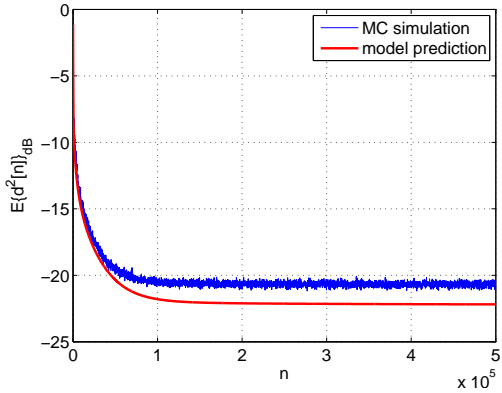}
}
\caption{MC simulation for Model Verification 3}
\label{fig:nonstat}
\end{figure}

The model in this paper is derived under stationarity assumptions for the input signals. Nevertheless, for nonstationary signals, it preserves sufficient information about the adaptation process to derive useful design guidelines. It is important to stress the fundamental difference between design guidelines and design rules.
Design guidelines are not rules that should be followed to design the system with a desired exact performance. 
No stochastic model can provide such rules as  
analytical models for adaptive algorithm behavior always rely on assumptions needed for mathematical tractability.  
Though the stationarity assumption is not satisfied in most practical systems, it is largely recognized that the models derived using them can still show tendencies of the algorithm behavior for reasonably small degrees of nonstationarity.
To illustrate the validity of the model even for a nonstationary input signal, $3$ simulation scenarios were tested for a GSC-AEC system with $\M=2$, $\Nh=1024$, $\NBF=16$, $\NAEC=1039$, $\muAEC=\muBF=3.1612\times 10^{-4}$ and different degrees of nonstationarity $\eta$~\cite{Manolakis:2000:SAS}. 
Results are shown in Fig.~\ref{fig:nonstat}.
For higher the degrees of nonstationarity, the behavior of the system diverges from the theoretical prediction. However, the system performance is still close enough to the MC simulation to jumpstart the design choices. 

\subsection{Design Example 1}
\label{sec:examples}

Consider an acoustic echo cancellation system with a reverberation time $T_{\rm R}(60) = 60$ ms (typical of a car cabin) and a background noise level of $-20$ dB. We assume the frequency response in the desired DOA (broadside) can be modeled by a linear phase, delayless all-pass filter with $\NfCal=16$ coefficients, i.e. $\fCal = [1,\mO_{1 \times \NfCal -1}]^T$~\cite{Frost:1972:AALCAAP}. The design goals are convergence of $J[n]_{\mathrm{dB}} < -20$ in less than $2$~s (evaluated at $n=1.5E4$ with \eqref{eqn:rhorecursivKailath} and \eqref{eqn:Jmsweighterror}).

We consider the design using the step-size matrix in~\eqref{eqn:stepsizematrix}. 
The frequency response model requires $\NBF=\NfCal=16$. The LEM plants have length $\Nh > f_{\rm s} T_R(60)=480$, and we thus set $\Nh=500$.
The free parameters are then $\M$ and $\NAEC$. The choices of $\M$ and $\NAEC$ affect the computational complexity per iteration and the convergence speed.

\begin{figure}[ht!]
\centering
   \includegraphics[width = 0.45\textwidth]{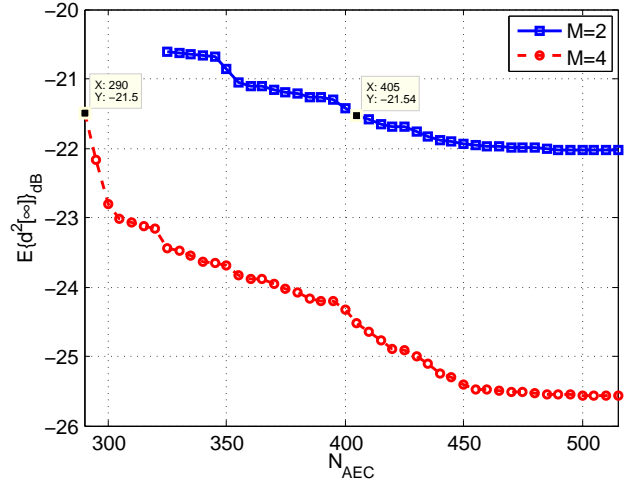}
 \caption{$E\{d^2[\infty]\}$ versus $\NAEC$ for different values of $\M$. Input is AR1(-0.9).}
\label{fig:perfcurve}
\end{figure}

Fig.~\ref{fig:perfcurve} was produced evaluating~\eqref{eqn:rhorecursivKailath} and~\eqref{eqn:Jmsweighterror} 
for $\M=\{1,2\}$, $\NAEC$ from $290$ to $515$, and $\ftr(\mM \mRmod)= \{2/300, \ldots, 2/3\}$ to compute $J[\infty]$. If $J[\infty] > -20$ dB this configuration is discarded.
 ~\eqref{eqn:rhorecursivKailath} and~\eqref{eqn:Jmsweighterror} are evaluated for $\muAEC = \{ 0.01, \ldots 0.99 \} \times 2(J[\infty] - \Jmin)/(J[\infty] \ftr(\mRuhcuhc))$.
If more than one combination of $(\muAEC,\muBF)$ is capable of reaching the desired cancellation at $n=1.5E4$ then only the one with the lower $J[\infty]$ is considered. 

\begin{figure}[ht!]
\centering
   \includegraphics[width = 0.45\textwidth]{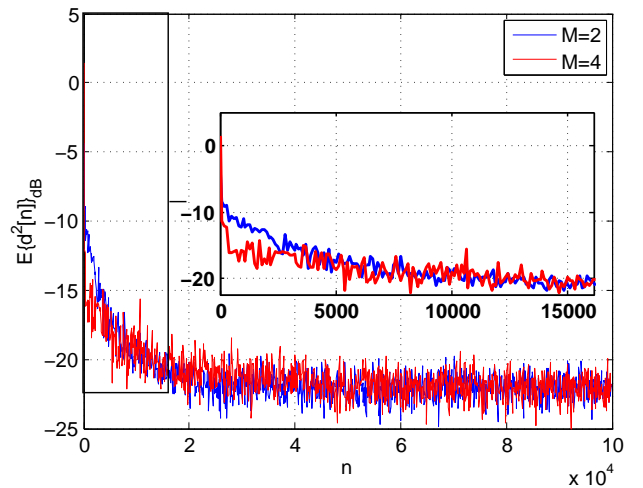}
 \caption{MC simulation (speech signals)}
\label{fig:perfcurve2}
\end{figure}

From Fig.~\ref{fig:perfcurve}, 2 candidate solutions were selected with $J[\infty] \approx -21.5$ dB and simulated using real speech signals (with pauses removed). 
Average results for an ensemble of $50$ runs are shown in Fig.~\ref{fig:perfcurve2}. 
\begin{table}[ht!]
\center
\caption{Parameters for Design Example 1}
\begin{tabular}{rrrrr}
\hline
\multicolumn{1}{l}{$M$} & \multicolumn{1}{l}{$\NAEC$} & \multicolumn{1}{l}{$\muAEC$} & \multicolumn{1}{l}{$\muBF$} & \multicolumn{1}{l}{$J[\infty]$} \\\hline
$2$ & $405$ & $9.7778E-04$ & $6.4603E-04$ & $-21.54$ dB\\
$4$ & $290$ & $9.1034E-04$ & $1.6969E-04$ & $-21.5$ dB\\\hline
\end{tabular}
\end{table}

\subsection{Design Example 2}

\begin{figure}[ht!]
\centering
   \includegraphics[width = 0.45\textwidth]{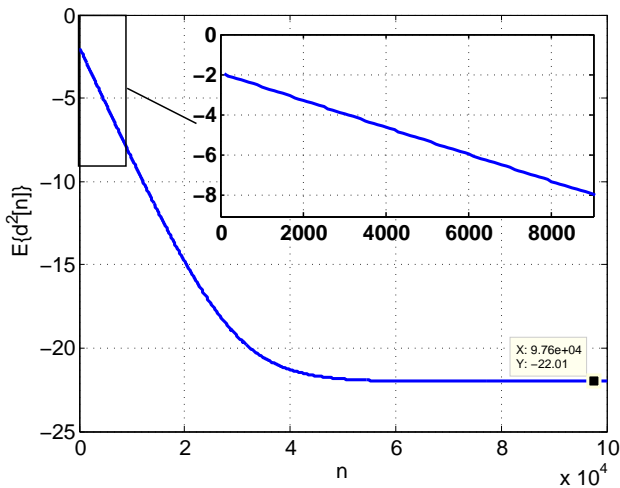}
 \caption{Model prediction for Design Example 2}
\label{fig:whitening}
\end{figure}

To illustrate the use of the derived optimal step-matrix derived in~\eqref{eqn:optstepsizemat}, consider an acoustic echo canceler with $\M=2$ microphones for a large conference room with reverberation time $T_{\mathrm{R}}(60)=0.45$ s. Assuming a sampling rate of $8$ kHz, the LEM plant length is about $\Nh \approx 0.45. 8 \times 10^3 = 3600$ coefficients. For this sampling rate, we also consider the frequency response in the desired DOA can be guaranteed with $\NfCal = 16$ constraints. In this simulation we consider only design choices with $\NAEC=\Nh+\NBF-1$.
The far-end signal is modeled by a unity variance AR1($-0.9$) random process and the noises in each microphone are assumed independent and  modeled by Gaussian i.i.d. variables with variance $10^{-2}$. In this design we desire a steady-state MOP of $-22$ dB and a $-10$ dB MOP after  $1$ second of convergence ($n=8000$).  To verify the feasibility of this design we used the proposed statistical model to predict the transient behavior using the optimal step-matrix derived in~\eqref{eqn:optstepsizemat}. Results are shown in Fig.~\ref{fig:whitening}. 
We observe that the whitening transformation is capable of achieving a MOP around $-8$ dB at $n=8000$. As this is the design with the optimal convergence speed, it is not possible to achieve, at the same time, the desired steady-state MOP of $-22$ dB and less than $-8$ dB at $n=8000$. Therefore it is not possible to design an BF-AEC system  in the GSC form with $T_{\mathrm{R}}(60)=0.45$ and these performance requirements.


\subsection{Design Example 3}
\label{subsec:controlstates}

\begin{figure}[htb]
  \centering
\includegraphics[width=8cm]{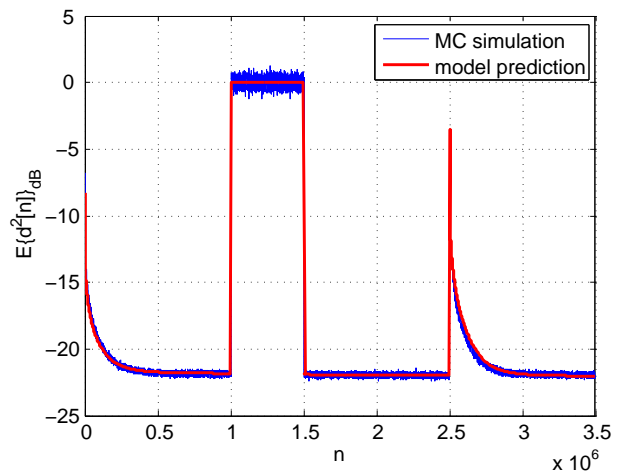}
\caption{Model prediction and Monte-Carlo simulation results Design Example 3}
\label{fig:controlstates}
\end{figure}

Consider an BF-AEC system designed to work in a room with reverberation time $T_R(60) = 100$ ms. To guarantee a Public Switched Telephone Network (PSTN) quality signal, the sampling rate was chosen $f_{\mathrm{s}}=8000$. 
To model the LEM plant impulse responses we have chosen $\Nh=1000$ coefficients.
Assume there are $\M=2$ microphones available and a reasonable frequency response in the DOA is achieved with $\NBF=16$. Finally, an AEC length $\NAEC=\Nh + \NBF -1 = 1015$ coefficients was used. 
For this example, we assume an ideal double-talk detector is available.
The initial DOA is assumed initially at a $\pi/4$ angle in relation to the broadside of the microphone array. During the first $10^6$ samples, the adaptation occurs in the absence of near-end speech with equal step-sizes $\muAEC = \muBF = 8.2147 \times 10^{-5}$. 
Then, a double-talk period occurs. 
The unitary power near-end speech, modelled as an AR1(-0.9) process, arrives from the broadside of the microphone array. 
We assume double-talk control logic, constraint and block matrix correction act instantly. 
The adaptation of the AEC is frozen by the double-talk detector ($\muAEC=0$), and during the next $5.10^{5}$ samples only the BF is adapted with $\muBF = 8.2147 \times 10^{-5}$. 
During this period, the convergence is significantly faster as the BF-AEC structure is not jointly-optimized and the effective adaptive filter length is reduced to $\M \NBF - \NfCal$.
During the next $10^6$ samples, the near-end speech is removed and, considering the BF state from the double-talk period is a good initial solution, step-sizes are set to 
$\muAEC = 9.3243\times 10^{-5}$ and $\muBF= 10^{-7}$. Finally, the LEM plant is subjected to an abrupt change in which a completely new $\mH$ is used. In this configuration the new step-sizes are set to $\muAEC=\muBF = 9\times 10^{-5}$ to accelerate convergence after an abrupt LEM plant change~\cite{Tourneret:2009:ECTGLRTDTCC}. The model predictions and Monte Carlo Simulations (ensemble of 50 runs) are compared on Fig.~\ref{fig:controlstates}

\section{Conclusion}

This work presented a statistical analysis of a class of jointly optimized beamformer-assisted AEC.  
The analysis was performed for systems with the BF implemented in the GSC form and using the LMS algorithm. The analysis considered convergence control using a step size matrix to accommodate typical control logic implementations.  We have shown that the joint optimization of the BF-AEC  is equivalent to a LCMV problem. Thus, the derived analytical models can be used to predict the transient performance of general adaptive wideband beamformers. The stochastic model was determined for the transient and steady-state behaviors of the residual mean echo power for stationary Gaussian inputs. Convergence analysis lead to stability bounds for the step-size matrix.  Design guidelines were derived from the analytical models. Monte Carlo simulations illustrated the accuracy of the theoretical models and the applicability of the proposed design guidelines.  
Finally, it was shown how a high convergence rate can be achieved using a quasi-Newton  adaptation scheme in which the step-size matrix is designed to whiten the combined input vector.

\section*{Acknowledgements}
 
 The authors would like to thank the invaluable help of Prof. Abraham Alcaim, from CETUC at PUC-Rio, who generously provided the speech signal database from which the real speech signals used in the described experiments were taken.

\bibliographystyle{IEEEtran}
\bibliography{acousticecho,microphonearrays}

\end{document}